\documentclass[a4paper,11pt, twoside, leqn, draft]{article}

\usepackage{amsmath}
\usepackage{amsfonts}
\usepackage{amssymb}
\usepackage{amsthm}
\usepackage{stackrel}
\usepackage[nobysame]{amsrefs}
\usepackage[cp1250]{inputenc}
\usepackage{setspace}
\usepackage[T1]{fontenc}
\usepackage{geometry}
\usepackage{fancyhdr}
\usepackage{graphicx}
\usepackage{color}
\usepackage{enumerate}
\usepackage{fullpage}
\newcommand{\exe}{\mathcal{E}}
\newcommand{\Te}{B_{q'}}
\newcommand{\eN}{{\ensuremath{\mathbb N}}}
\newcommand{\Ex}{{\ensuremath{\mathbb E}}}
\newcommand{\Pro}{{\ensuremath{\mathbb P}}}
\newcommand{\eR}{{\ensuremath{\mathbb R}}}
\newcommand{\lv}{\left\lVert}
\newcommand{\rv}{\right\rVert}
\newcommand{\1}{\mathbb{1}}
\newcommand{\wy}{\mathcal{E}}

\newcommand{\eps}{\varepsilon}

\newcommand{\te}{\otimes}
\newcommand{\bb}{\beta_{A,S}}

\newcommand{\ve}{\varepsilon}
\newcommand{\norxy}[1]{\lv #1 \rv_{X,Y,p}}
\newcommand{\su}{\sup_{t\in \Te}}
\newcommand{\mi}{\bar{\mu}_{n,\eps}}

\newcommand{\supf}{\sup_{f\in B^*(L_q)}}
\newcommand{\supt}{\sup_{f\in B^*(F)}}
\newcommand{\ban}{(F,\lv \cdot \rv)}
\newcommand{\kul}{B^*(L_q)}
\newcommand{\kxx}{B^X_p}
\newcommand{\kyy}{B^Y_p}
\newcommand{\wspol}{s_A(\Te)}
\newcommand{\alinf}{\alpha_{\infty,A}}
\newcommand{\dinf}{d_{\infty,A}}
\newcommand{\diaminf}{\Delta_{\infty,A}}

\newcommand{\alf}{\alpha_{A}}

\newcommand{\beet}{\beta_{A,\Te}}

\newcommand{\de}{\Delta} \newcommand{\hd}{\tilde{\Delta}}

\newcommand{\nx}{\hat{N}^X_i}
\newcommand{\my}{\hat{N}^Y_j}
\newcommand{\kx}{\sum_{i} \nx(x_i) \leq p}
\newcommand{\ky}{\sum_{j} \my(y_j) \leq p}
\newcommand{\bt}{\beta_{A,S}}
\newcommand{\btt}{\beta_{A,T}}

\newcommand{\da}{\Delta_A}
\newcommand{\ten}{\otimes}
\newcommand{\nory}[1]{\lv #1 \rv_{Y,p}}
\newcommand{\norx}[1] {\lv #1 \rv_{X,p}}
\newtheorem{theorem}[subsection]{Theorem}
\newtheorem{lemma}[subsection]{Lemma}
\newtheorem{proposition}[subsection]{Proposition}

\newtheorem{fact}[subsection]{Fact}
\newtheorem{corollary}[subsection]{Corollary}
\newtheorem{remark}[subsection]{Remark}
\providecommand{\keywords}[1]{\textbf{\textit{Keywords: }} #1}

\providecommand{\klas}[1]{\textbf{\textit{AMS MSC 2010: }} #1}

\begin{document}
\author{Rafa{\l} Meller\footnote{Institute of Mathematics, University of Warsaw, Banacha 2, 02-097 Warsaw, Poland, 
    rmeller@mimuw.edu.pl}\\
Institute of Mathematics\\
University of Warsaw\\
02-097 Warszawa, Poland\\
E-mail: r.meller@mimuw.edu.pl}
\title{Moments and tails of $L_q$-valued chaoses of order two based on independent variables with log-concave tails\thanks{This research was partly supported by the National Science Centre, Poland grants 2018/29/N/ST1/00454 and  2021/40/C/ST1/00330}}
\date{}
\maketitle

\renewcommand{\thefootnote}{}



\renewcommand{\thefootnote}{\arabic{footnote}}
\setcounter{footnote}{0}


\begin{abstract}
We derive a lower bound for moments of random chaoses of order two with coefficients in arbitrary Banach space $F$ generated by independent symmetric random variables with logarithmically concave tails (which is probably two-sided). We also provide two upper bounds for moments of such chaoses when $F=L_q$. The first is valid under the additional subgaussanity assumption. The second one does not require additional assumptions but is not optimal in general. Both upper bounds are sufficient for obtaining two-sided moment estimates for chaoses with values in $L_q$  generated by Weibull random variables with shape parameter greater or equal to $1$.

\keywords{Random chaoses; Random quadratic forms; Tail and moment estimates; Logarithmically concave tails} \\
\klas{60E15}
\end{abstract}

\section{Introduction}

A (homogeneous) polynomial chaos of order $d$ is a random variable defined as
\[S=\sum_{i_1,\ldots,i_d} a_{i_1,\ldots,i_d} X_{i_1}\cdots X_{i_d},\]
where $X_1,\ldots,X_n$ are independent random variables and $(a_{i_1,\ldots,i_d})_{1\leq i_1,\ldots,i_d\leq n}$ are  coefficients ($d$-indexed) that belong to a Banach space $(F,\lv \cdot \rv)$ such that $\ a_{i_1,\ldots,i_d}=0$  if  $i_k=i_l$ for some  $k\neq l$. The natural question is whether it is possible to give an exact description of the moments of $S$ defined as $\lv S \rv_p:=(\Ex \lv S \rv ^p)^{1/p}$. By exact bounds, we mean two-sided ones. So, we look, for deterministic expression $f$ (which  may depend on $p$, the coefficients etc.) such that
\begin{equation}\label{prob1}
\frac{f}{C(d)}\leq \lv S \rv_p \leq C(d)f,   
\end{equation}
where $C(d)$ is some constant that depends only on $d$, the order of chaos. To derive an effective formula for $f$ one must assume something about the distribution of $(X_i)_i$ (with an exception on the real line for $d=1$). For several reasons, it is very convenient to work with symmetric random variables $X_i$'s that have logarithmically concave tails (LCT for short) i.e., such that for any $i$ the function $t\mapsto -\ln \Pro(|X_i|\geq t) \in [0,\infty]$ is convex. First, this class contains many natural distributions such as Gaussian, exponential, Rademacher, and Weibull (with a shape parameter greater than or equal to 1). Second, in this class \eqref{prob1} can be proved with a reasonable formula $f$ in several cases (see examples below). Finally, in this class, moments estimates imply bounds on the tails of $S$ that is on $\Pro(\lv S \rv \geq t),\ t>0$. The argument is standard regardless of the chaos order and the Banach space's choice.\\
The problem of establishing \eqref{prob1} naturally falls into two parts depending on whether the coefficients $a_{i_1,\ldots,i_d}$ are real or vector-valued,  with significant unanswered questions even on the real line.\\
Now we present the state of the art. We start with the real coefficients. For $d=1$, deterministic bounds on the moments of $S$ are known only under the symmetry or non-negativity assumption \cite{d1}.  If the variables and coefficients are non-negative, then the formula for the moments is known for any order $d$ under the assumption that the moments of the r.v.'s $X_i$ grow at most polynomially (moments of a variable with logarithmically concave tails grow at most linearly) \cite{ja}. The case of symmetric variables is much less understood. We know formulas for $d=2$ when the moments of variables grow at most polynomially \cite{MR3984282}, and for $d=3$, when the variables have logarithmic tails \cite{adlat}. For any $d$, formulas are known if the variables are Gaussian \cite{latgaus}, exponential \cite{adlat}, or have logarithmically convex tails \cite{kol} (which includes the Weibull random variable with a shape parameter less than or equal to $1$). In the Gaussian case, something even more general is known, namely the formula for the moments of any polynomial in Gaussian variables \cite{adwolf}. {\L}ochowski's result \cite{loch} is also worth mentioning. He obtained bounds for chaos of arbitrary order based on symmetric r.v.'s with logarithmically concave tails. However, his bounds involve expectations for the suprema of non-Gaussian processes, which are very difficult to estimate.  \\
Very little is known in the case of vector-valued coefficients, except for chaos of order $1$. If $d=1$, we know formulas when the moments of $(X_i)_i$ grow at most polynomially \cite{strzelec}, i.e. in quite satisfying generality. For $d\geq 2$, nothing is known except for two cases: Gaussian chaos of arbitrary order in spaces of type $2$ and exponential chaos of arbitrary order in $L_q$  spaces \cite{adlatmel,gaban}. The latter results could only be derived because the exponential symmetric r.v. is (almost) a product of two independent Gaussian r.v.'s. Thus, we could reduce the study of exponential chaos of order $d$ to Gaussian chaos of order $2d$. There are formulas in the literature for moments of Gaussian chaos of any order in arbitrary Banach space. Unfortunately, they involve the expectation of suprema of Gaussian processes (which are very hard to estimate), see \cite[Chapter 3.2]{proinban} and \cite{Borell} for details. It is also worth mentioning Adamczak's result \cite{Adban}. Although his bounds are not two-sided (which was not the goal of this work), they are obtained in greater generality and any Banach space. To some extent, our formulas are similar to those obtained by Adamczak.\\

We are interested in studying moments of $L_q$-valued random chaos of order 2, based on symmetric random variables with logarithmically concave tails. This can be seen as the first attempt to study the moments in the non-Gaussian case (as we mentioned, the exponential chaoses can be reduced to the Gaussian ones). In our setup
\[S=\sum_{i,j=1}^n a_{ij} X_i X_j,\ \textrm{where } a_{ij} \in L_q,\textrm{ and }  a_{ii}\equiv 0, \]
where $X_1,X_2,\ldots$ are independent symmetric random variables with LCT (with logarithmically concave tails). If $a_{ii} \neq 0$ then standard arguments like Fact \ref{fak:przek} allow us to study the diagonal and off-diagonal parts of the chaos separately. The diagonal part is a linear combination of independent random variables, and certain results can be applied, e.g. \cite{strzelec}. Thus, we focus on the case when $a_{ii}\equiv 0$. Without loss of generality, we can assume that $a_{ij}=a_{ji}$ and then use Theorem \ref{twr:dec} to reduce the problem of estimating $\lv S \rv_p$ to estimating moments of decoupled chaos
\[S'=\sum_{i,j=1}^n a_{ij} X_i Y_j,\ a_{ii}\equiv 0, \]
where the random vector $(Y_i)_i$ is an independent copy of the random vector $(X_i)_i$ (\cite{dec1,dec2}). The latter object has a richer structure and is easier to work with. Therefore, in this paper, we will focus only on decoupled chaos. We will also not assume that $a_{ii}\equiv 0$ since this is irrelevant in the case of decoupled chaos.\\
First, we give a simple lower bound for $\lv S'\rv_p$ in any Banach space $\ban$, which is probably two-sided (at least in a class of Banach spaces with non-trivial cotype). Then we restrict our attention to $L_q$ spaces. We give two upper bounds. The first is optimal in the class of sub-Gaussian random variables with logarithmically concave tails. The second is not optimal, but it does not require additional assumptions on the random variables. These theorems imply two-sided moment estimates for $L_q$ valued chaos based on symmetric random variables with CDF equal to $1-\exp(-x^r)$ for $r\geq 1$ (Weibull random variable with shape parameter $r$). This can be considered as the main result of the paper.

In the next section, we set up a notation and present the main results. In Section \ref{r5s2} we show that the main difficulty in obtaining upper bounds for the moments of random chaos is to properly estimate the expectation of the supremum of a given stochastic process. In Section \ref{rgaus} we derive the upper bound for the expectation of the supremum of a certain Gaussian process. It is a generalization of the upper bound obtained in \cite{adlatmel}. Section \ref{rexp} deals with the case of the supremum of an exponential process. In Section \ref{rrozbij} we use the ideas from \cite{adlat} to derive the decomposition theorem for exponential processes in $L_q$ space. In the last section, we prove the main theorems.  Unfortunately, our approach is not sufficient to obtain a more general result, and new ideas are needed. In this paper, we use many results from other papers. To overcome this inconvenience, we have collected all the cited facts in the appendix. We also provide a glossary explaining the notation that precedes the bibliography.

\textbf{Acknowledgments}
 I want to thank prof. R. Lata{\l{}}a for pointing out that Theorems \ref{r5thmmain}, \ref{Thmnewres} imply Theorem \ref{thmwyk}.

\section{Notation, convention and main results} \label{r5s1}

In this note, $g_1,g_2,\ldots,\wy_1,\wy_2,\ldots$ denote independent random variables with standard Gaussian and symmetric exponential distributions (i.e. the distribution with density  $1/2 \exp(-|x|/2)$). Here and subsequently $G_n$ stands for $(g_1,\ldots,g_n)$ and $E_n$ for $(\wy_1,\ldots,\wy_n)$.
By $\Ex^X,\Ex^Y$ we mean integration with respect to $X_1,X_2,\ldots$ and $Y_1,Y_2,\ldots$ respectively.

The letter $C$ (resp. $C(\alpha)$) stands for a numerical constant (resp. constant that depends only on some parameter $\alpha$), which may be different at each occurrence. We  use the notation $a \lesssim b$ (resp. $a \lesssim^\alpha b$) if $a \leq Cb$ (resp. $a\leq C(\alpha) b$). We will also write $a \approx b$ (resp. $a \approx ^\alpha b$) if $a\lesssim b$ and $b\lesssim a$ (resp. $a \lesssim^\alpha b$ and $b \lesssim^\alpha a$).  If $(F,\lv \cdot \rv)$ is a Banach space, then $(F^*,\lv \cdot \rv_*)$ stands for its dual space and $B^*(F)=\{ f\in F^*: \lv f \rv_*\leq 1\}$ for the unit ball in the dual space. If $q
\geq 1$, then by $q'$ we denote the H\"{o}lder conjugate number to $q$ (i.e. $1/q+1/q'=1$).
We also denote $B^n_{u}=\{v\in \eR^n: \sum_i^n |v_i|^{u}\leq 1\}$. We often omit the upper index $n$ if it does not lead to misunderstandings.
Obviously, $B^n_{q'}=B^*(\ell_q^n)$, where $\ell_q^n$ is the space of sequences of length $n$ equipped with the norm $\lv x \rv^q_{\ell_q}=\sum_i^n |x_i|^q$.\\
Let $(X_i)_i,(Y_j)_j$ be independent, symmetric random variables. We define the functions
\[N^X_i(t)=-\ln \Pro(|X_i|\geq t) \in [0,\infty],\ N^Y_j(t)=-\ln \Pro(|Y_j|\geq t)\in [0,\infty].\]
Our basic assumption is that $(X_i)_i,(Y_j)_j$ have LCT (logarithmically concave tails) i.e. the functions $(N^X_i)_i,(N^Y_j)_j$ are convex. We are mainly concerned with the homogeneous inequalities of degree one, so we can normalize the random variables as follows
\begin{align} \label{normalizacja}
\inf \left\{ t\geq 0,\ N^X_i(t)\geq 1\right\}=\inf \left\{ t\geq 0,\ N^Y_i(t)\geq 1\right\}=1.
\end{align}

\noindent We set
$$\nx(t)=\begin{cases}t^2 &\textrm{for } |t|\leq 1 \\ N^X_i(t) &\textrm{for } |t|>1 \end{cases},\quad \my(t)=\begin{cases}t^2 &\textrm{for } |t|\leq 1 \\ N^Y_i(t) &\textrm{for } |t|>1 \end{cases}.  $$

\noindent Note that the convexity of $N^X_i, N^Y_j$ and the normalization condition \eqref{normalizacja} imply that
\begin{equation}
\nx(t)=N^X_i(t)\geq t,\quad \my(t)=N^Y_i(t)\geq t\quad \textrm{for } t\geq 1,\label{r5cos}
\end{equation}
\begin{equation}
1/e \leq \Ex X^2_i, \Ex Y^2_j \leq 1+4/e\leq 3,\quad \Ex X^4_i,\Ex Y^4_j \leq 1+64/e. \label{wnnorm}
\end{equation}
The first formula is clear, to prove the second one it is enough to observe that 
$$1/e \leq \Pro(|X^2_i|\geq 1) \leq \Ex X^2_i\leq 1+\int_{1}^\infty 2xe^{-x} dx =1+4/e\leq 3$$
(we prove the bounds for $\Ex X^4_i$ analogously).
We define
\begin{equation}\label{kulax}
\kxx=\{x\in \eR^n: \sum_i \nx(x_i) \leq p \},\quad \kyy= \{y \in \eR^n: \sum_j \my(y_j) \leq p\}.
\end{equation}
Observe that \eqref{r5cos} implies that
\begin{equation}\label{r5zawieranie}
\kxx,\kyy \subset   p^{1/2}  B^n_2+pB^n_1.
\end{equation}
Let $(a_{ij})_{ij}$ be an $\eR$ valued matrix and $(a_i)_i \in \eR^n$.
 The following three norms will play a crucial role in this paper:
\begin{align*}
\lv (a_{ij})_{ij} \rv_{X,Y,p}&:=\sup \left\{  \sum_{ij} a_{ij} x_i y_j   : \kx,\ \ky \right\}=\sup_{\substack{x \in \kxx\\y\in \kyy}} \sum_{ij} a_{ij} x_iy_j, \\
\lv (a_i)_i \rv_{X,p} &:= \sup \left\{  \sum_i a_i x_i   : \kx \right\}=\sup_{x \in \kxx} \sum_i a_i x_i ,\\
\lv (a_j)_j \rv_{Y,p} &:= \sup \left\{ \sum_j a_j y_j \:  \ky \right\}=\sup_{y \in \kyy} \sum_j a_j y_j. 
\end{align*}
It can be shown that the above objects are norms, not just seminorms (this follows from our normalization \eqref{r5cos}). \\
Since $\nx(t/u)\leq \nx(t)/u$ for $u\geq 1$ (it follows from the convexity of $N^X_i$, the normalization condition \eqref{normalizacja} and that $\nx(0)=0$), we have that
\begin{equation}\label{r5skalowanie}
\lv (a_i)_i \rv_{X,up} \leq u\norx{(a_i)_i},\   \lv (a_{ij})_{ij} \rv_{X,Y,up}\leq u^2\norxy{(a_{ij})_{ij}},
\end{equation}
and the first formula is also valid for $\nory{(a_j)_j}$.\\
\noindent We start with a simple lower bound.
\begin{proposition}\label{r5zdolu}
Assume that $X_1,X_2,\ldots,Y_1,Y_2,\ldots$ are independent symmetric r.v.'s with LCT such that the normalization condition \eqref{normalizacja} holds. Let $(a_{ij})_{ij}$ belong to  a Banach space $\ban$. Then for any $p\geq 1$ we have
\begin{align}
&\lv \sum_{ij} a_{ij} X_i Y_j \rv_p \gtrsim  \Ex\lv \sum_{ij} a_{ij} X_i Y_j \rv  + \sup_{x\in \kxx}\Ex \lv\sum_{ij} a_{ij} x_i Y_j \rv +\sup_{y \in \kyy} \Ex \lv \sum_{ij} a_{ij} X_i y_j \rv
  \nonumber \\
&+\supt \norx{\left( \sqrt{\sum_j (f(a_{ij}))^2}\right)_i}+\supt \nory{\left( \sqrt{\sum_i (f(a_{ij}))^2}\right)_j}\nonumber \\
&+\supt\norxy{( f(a_{ij}))_{ij} }=:W^{X,Y}(p).\label{r5indol}
\end{align}
Moreover, for any $t\geq 1$
\begin{align}\label{ogonbound}
\Pro\left(\lv \sum_{ij} a_{ij} X_i Y_j \rv\geq W^{X,Y}(t)\right) \gtrsim \exp(-Ct).
\end{align}
\end{proposition}
\begin{proof}
Trivially $\lv \sum_{ij} a_{ij} X_i Y_j  \rv_p \geq \Ex \lv \sum_{ij} a_{ij} X_i Y_j \rv$. Using Theorem \ref{r5some1} we obtain
\begin{align*}
\lv  \sum_{ij} a_{ij} X_i Y_j  \rv_p &=\lv \supt \sum_{ij} f(a_{ij})X_iY_j \rv_p  \geq \supt \lv \sum_{ij} f(a_{ij}) X_i Y_j \rv_p \\
&\gtrsim \supt \norx{\left(\sqrt{\sum_j(f( a_{ij}))^2} \right)_i}+\supt \nory{\left( \sqrt{\sum_i (f(a_{ij}))^2}\right)_j}\\
&+\supt\norxy{\left( f ( a_{ij})\right)_{ij}}.
\end{align*}
Jensen's inequality and Theorem \ref{r5gl}  yield
\begin{align*}
&\lv  \sum_{ij} a_{ij} X_i Y_j  \rv_p  \geq \left( \Ex^Y \supt \Ex^X \left|\sum_{ij} f(a_{ij}) X_i Y_j  \right|^p \right)^{1/p} \\
&\gtrsim \lv \sup_{f\in B^*(F), x \in \kxx} \sum_{ij} f(a_{ij}) x_i Y_j  \rv_p \geq \sup_{x \in \kxx} \lv \sum_{ij} a_{ij} x_i Y_j\rv_p \geq \sup_{x \in \kxx} \Ex \lv \sum_{ij} a_{ij} x_i Y_j\rv.
\end{align*}
Analogously $\lv  \sum_{ij} a_{ij} X_i Y_j  \rv_p \gtrsim\sup_{y \in \kyy} \Ex \lv \sum_{ij} a_{ij} X_i y_j \rv$ and  \eqref{ogonbound} follows. A standard application of the Paley-Zygmund inequality, the fact that the $p$th and $2p$th moments of $\sum_{ij} a_{ij} X_i Y_j$ are comparable up to a numerical constant (by Fact \ref{hiper}) and \eqref{r5skalowanie} imply \eqref{ogonbound} (for details see the proof of \cite[Corollary 2.3]{MR3984282}).
\end{proof}

\noindent  We suspect that the inequality \eqref{r5indol} can be reversed when the Banach space has the non-trivial co-type. Whether it can be reversed in any Banach space would be a challenging problem. A positive premise is that \eqref{r5indol} can be reversed in any Banach space with logarithmic accuracy when the variables have the normal distribution cf. \cite[Theorem 3]{adlatmel}. However, there may also be a simple counterexample. In this paper, we provide some results for $L_q$ spaces. In particular, we rely heavily on the fact that these spaces have a non-trivial co-type.\\
First, we formulate an upper bound for the moments of chaos under an additional assumption of subgaussianity. We say that a r.v. $Y$ is subgaussian with constant $\gamma>0$ if $\Ex Y=0$ and for any $t\in \eR$,  $\Ex \exp(tY) \leq \exp(\gamma t^2)$. 
\begin{theorem}\label{r5thmmain}
Assume that $X_1,X_2,\ldots,Y_1,Y_2,\ldots$ are independent symmetric r.v.'s with LCT such that \eqref{normalizacja} holds. Assume also that  $Y_1,Y_2,\ldots$ are subgaussian with constant $\gamma$. Let $(a_{ij})_{ij}$ be an $L_q(V,\mu)$-valued matrix. Then for any $p\geq 1$ we have
\begin{align}
&\lv  \sum_{ij} a_{ij} X_i Y_j   \rv_p \lesssim^q  \gamma \Ex\lv \sum_{ij} a_{ij} X_i Y_j \rv_{L_q}  + \sup_{\substack{x\in \kxx\\y \in \kyy} }\left(\Ex \lv\sum_{ij} a_{ij} x_i Y_j \rv_{L_q} +\Ex \lv \sum_{ij} a_{ij} X_i y_j \rv_{L_q}\right)  
 \nonumber \\
&+\supf \norx{\left( \sqrt{\sum_j (f(a_{ij}))^2}\right)_i}+\supf\norxy{( f(a_{ij}))_{ij} }=:W^{X,Y,\gamma}(p),\label{r5:main:ineq}
\end{align}
and
\begin{equation}\label{ogonbound2}
\Pro\left(\lv \sum_{ij} a_{ij} X_i Y_j \rv\geq W^{X,Y,\gamma}(t)\right) \leq C(q)\exp\left(-\frac{W^{X,Y,\gamma}(t)}{C(q)}\right).
    \end{equation}
\end{theorem}
\begin{remark}
A version of Lemma \ref{turbolemat} (which is for $\ell_q$ spaces, but the proof for $L_q$ spaces is identical) implies that
\begin{align*}
&\Ex\lv \sum_{ij} a_{ij} X_i Y_j \rv_{L_q} \approx^q \lv \sqrt{\sum_{ij} a^2_{ij} } \rv_{L_q},\ \Ex \lv\sum_{ij} a_{ij} x_i Y_j \rv_{L_q} \approx^q  \lv \sqrt{\sum_{j} \left(\sum_i x_i a_{ij}\right)^2 } \rv_{L_q},\\
&\Ex \lv\sum_{ij} a_{ij} X_i y_j \rv_{L_q} \approx^q  \lv \sqrt{\sum_{i} \left(\sum_i y_j a_{ij}\right)^2 } \rv_{L_q}.
\end{align*}
Thus one may rewrite \eqref{r5:main:ineq} without using expectations.  
\end{remark}
\noindent The formula \eqref{ogonbound2} is a simple consequence of Chebyshev's inequality and  \eqref{r5:main:ineq}. We will prove the latter in Section \ref{r5s4}. Proposition \ref{r5zdolu} ensures that \eqref{r5:main:ineq} and \eqref{ogonbound2} are in fact two-sided. \\
We derive a similar result without assuming the subgaussianity at the cost of an extra term. The theorem below also implies a bound on tails in a standard
way.
\begin{theorem}\label{Thmnewres}
Assume that $X_1,X_2,\ldots,Y_1,Y_2,\ldots$ are independent symmetric r.v.'s with LCT such that \eqref{normalizacja} holds. Let $(a_{ij})_{ij}$ be an $L_q(V,\mu)$-valued matrix. Then for any $p\geq 1$ we have
\begin{align}
\lv  \sum_{ij} a_{ij} X_i Y_j  \rv_p \lesssim^{q} &\Ex \lv \sum_{ij} a_{ij} X_i Y_j \rv_{L_q}+\sup_{\substack{x\in \kxx\\y \in \kyy} }\left( \Ex \lv \sum_{ij} a_{ij} x_i Y_j \rv_{L_q}
+\Ex \lv \sum_{ij} a_{ij} X_i y_j \rv_{L_q}\right) \nonumber \\
&+\supf \norx{ \left(\sqrt{\sum_j f(a_{ij})^2}\right)_i}+\supf\norxy{(f(a_{ij}))_{ij}} \nonumber \\
&+p\max_i \supf \sqrt{\sum_j f^2(a_{ij})}.\label{niergor}
\end{align} 
\end{theorem}

\begin{remark}
If we exchange $(X_i)_i$ r.v.'s with $(Y_i)_i$ in Theorem \ref{Thmnewres}, we get a different upper bound. This is caused by the third and the last term in \eqref{niergor} which are not symmetric. The lack of symmetry regarding the third term is just a matter of formulation. In fact \eqref{niergor} can be shown with $\supf \norx{ \left(\sqrt{\sum_j f(a_{ij})^2}\right)_i}$ replaced by a smaller and symmetric quantity
\[\supf \inf_{I,J\subset \eN,|I|=|J|=p}\sqrt{\sum_{i \notin I, j \notin J} f(a_{ij})^2}.\]
But we decided to keep the presentation of \eqref{niergor} as it is. The second non-symmetric term is  $p\max_i \supf \sqrt{\sum_j f^2(a_{ij})}$. This is a remnant of our suboptimal proof. We suspect that a more subtle argument can eliminate this term.
\end{remark}

If $X_1,\ldots,Y_1,\ldots$ are independent Rademacher r.v.'s (symmetric $\pm 1$ r.v.'s) then 
 \[\lv  \sum_{ij} a_{ij} X_i Y_j  \rv_p \leq \lv  \sum_{ij} a_{ij} X_i Y_j  \rv_\infty <\infty,\]
 while
 \[p\max_i \supf \sqrt{\sum_j f^2(a_{ij})} \stackrel{p \to \infty}{\to} \infty.\]
So it is impossible to reverse the inequality \eqref{niergor}. However both Theorems \ref{r5thmmain} and \ref{Thmnewres} imply two-sided moments estimates for chaos based on Weibull r.v.'s with values in $L_q$ spaces. This can be considered as the main achievement of this paper. It is worth noting that the constants in the theorem below depend only on $q$ and not on $r$. Theorem \ref{thmwyk} implies a bound on the tails of the underlying chaos by standard arguments (cf. proof of Theorem \ref{r5zdolu} and comment after Theorem \ref{r5thmmain}).
\begin{theorem}\label{thmwyk}
Assume that $X_1,X_2,\ldots$ and $Y_1,Y_2,\ldots$ are independent symmetric random variables  with CDF given by $\Pro(|X_i|\geq t)=\Pro(|Y_i|\geq t)= \exp(-|x|^r)$, $r\geq 1$.  Let $(a_{ij})_{ij}$ be an $L_q(V,\mu)$-valued matrix. Then for any $p\geq 1$ we have
\begin{align}
&\lv  \sum_{ij} a_{ij} X_i Y_j \rv_p \approx^q \Ex\lv \sum_{ij} a_{ij} X_i Y_j \rv_{L_q}   + \sup_{\substack{x\in \kxx\\y \in \kyy} }\left(\Ex \lv\sum_{ij} a_{ij} x_i Y_j \rv_{L_q} +\Ex \lv \sum_{ij} a_{ij} X_i y_j \rv_{L_q}\right)
 \nonumber \\
&+\supf \norx{\left( \sqrt{\sum_j (f(a_{ij}))^2}\right)_i}+\supf\norxy{( f(a_{ij}))_{ij} }.\label{thmwyk:ineq}
\end{align}
\end{theorem}
\begin{proof}
The variables $X_1,X_2,\ldots,Y_1,Y_2,\ldots$ have LCT (since $r\geq 1$) and they satisfy \eqref{normalizacja}. So the lower bound follows from Proposition \ref{r5zdolu}. For $r\geq 2$ the variables in the theorem are subgaussian with constant $C=10$. So in this case we can use Theorem \ref{r5thmmain}. Now consider $r\in [1,2)$. Then $\nx(x)\leq x^2, \my(x)\leq x^2$ and as a result $p^{1/2}B_1\subset p^{1/2}B_2\subset \kxx,\ p^{1/2}B_2\subset \kyy$. Thus,
\begin{align*}
\supf \norxy{(f( a_{ij}))_{ij}}&\geq \sup\left\{ \sum_{ij} f(a_{ij})x_i y_j  :  x \in p^{1/2}B_1,\ y\in p^{1/2}B_2,\  f\in \kul \right\}.  
\end{align*}
By first taking supremum over $x$ and then over $y$ we get
\begin{align*}
\supf \norxy{(f (a_{ij}))}&\geq p^{1/2}\max_i \sup_{\lv y \rv_2=p^{1/2}, f\in \kul} \left| \sum_j f(a_{ij})y_j \right|\\
&= p\max_i \sup_{f\in \kul} \sqrt{\sum_j f^2(a_{ij})}.
\end{align*}
So the upper bound in \eqref{thmwyk:ineq} for  $r\in(1,2]$ follows by Theorem \ref{Thmnewres}.
\end{proof}

In Hilbert spaces, we do not need any additional assumptions.

\begin{theorem}\label{hil}
Assume that $X_1,X_2,\ldots$ and $Y_1,Y_2,\ldots$ are symmetric independent random variables with LCT such that \eqref{normalizacja} holds. Let $(a_{ij})_{ij}$ be an $(H,\lv \cdot \rv)$ -valued matrix, where $H$ is a Hilbert space. Then for any $p\geq 1$ we have
\begin{align}
&\lv   \sum_{ij} a_{ij} X_i Y_j   \rv_p \approx \Ex\lv \sum_{ij} a_{ij} X_i Y_j \rv  + \sup_{x\in \kxx}\Ex \lv\sum_{ij} a_{ij} x_i Y_j \rv 
 \nonumber \\
&+\sup_{y \in \kyy} \Ex \lv \sum_{ij} a_{ij} X_i y_j \rv+\sup_{f\in B^*(H)}\norx{\left( \sqrt{\sum_j (f(a_{ij}))^2}\right)_i}+\sup_{f\in B^*(H)}\norxy{( f(a_{ij}))_{ij} }.\label{r5hilnier}
\end{align}
\end{theorem}

\section{Reduction to a bound on the supremum of a certain stochastic process}
\label{r5s2}

We begin by showing that the only difficulty in bounding moments of $\lv\sum_{ij} a_{ij} X_i Y_j\rv$, is to bound the expectation of suprema of a certain stochastic process.
\begin{lemma}\label{r5redmom}
Assume that $X_1,X_2,\ldots$ and $Y_1,Y_2,\ldots$ are symmetric independent random variables with LCT such that \eqref{normalizacja} holds. Let $\ban$ be a Banach space and $(a_{ij})_{ij}$ be an $\ban$ -valued matrix. Then for any $p\geq 1$ we have
\begin{align*}
\lv  \sum_{ij} a_{ij} X_i Y_j \rv_p \lesssim& \ \Ex \lv \sum_{ij} a_{ij} X_i Y_j \rv +\sup_{y \in \kyy} \Ex \lv \sum_{ij} a_{ij} X_i y_j \rv +
\Ex \sup_{ x \in \kxx}  \lv \sum_{ij} a_{ij} x_i Y_j \rv  \\
&+ \supt \norxy{(f(a_{ij}) )_{ij}}.
\end{align*}
\end{lemma}
\begin{proof}
By conditionally applying Theorem \ref{r5lat1} we obtain
\begin{align*}
\lv  \sum_{ij} a_{ij} X_i Y_j \rv_p \lesssim \lv \left(\Ex^X \lv \sum_{ij} a_{ij} X_i Y_j \rv \right) \rv_p+ \lv \supt \norx{\left(\sum_{j}f(a_{ij})Y_j\right)_i} \rv_p.
\end{align*}
Since $y \mapsto  \Ex^X \lv \sum_{ij} a_{ij} X_i y_j\rv $ is a seminorm, Theorem \ref{r5lat1}  implies that
\begin{align*}
 \lv \left(\Ex^X \lv \sum_{ij} a_{ij} X_i Y_j \rv \right) \rv_p &\lesssim  \Ex \lv  \sum_{ij} a_{ij} X_i Y_j\rv  + \sup_{y \in \kyy} \Ex \lv\sum_{ij} a_{ij} X_i y_j \rv.
\end{align*}
We finish the proof by observing that $y\mapsto \supt \sup_{x\in \kxx}\sum_{ij} x_if(a_{ij})y_j$ is also a seminorm so again Theorem \ref{r5lat1} implies that
\begin{align*}
 \lv \supt \norx{ \left(\sum_{j} f(a_{ij})Y_j \right)_i} \rv_p&= \lv \supt \sup_{x\in \kxx}\sum_{ij} x_if(a_{ij})Y_j  \rv_p\\
&\lesssim \Ex \sup_{ x \in \kxx}  \lv \sum_{ij} a_{ij} x_i Y_j \rv  + \supt \norxy{(f(a_{ij}))_{ij}}.
\end{align*}
\end{proof}
So in order to prove  Theorems \ref{r5thmmain}, \ref{Thmnewres}, \ref{hil}   it is sufficient to establish upper bounds on $\Ex \sup_{x \in \kxx} \lv \sum_{ij} a_{ij} x_i Y_j \rv$. 
\begin{proof}[Proof of Theorem \ref{hil}]
 The lower bound follows by Proposition \ref{r5zdolu}.  It is sufficient to show the upper bound if we sum up over $i,j \leq n$. We can also assume that $\dim H<\infty$ so that w.l.o.g. $H=\eR^m$ and $\lv \cdot \rv$ is the Euclidean norm. Then $a_{ij}=(a_{ijk})_{k \leq m}$, $B^*(H)=B_2^m$ and $\lv x \rv=\sup_{t \in B_2} \sum_{i \leq m} t_i x_i$. Thanks to Lemma \ref{r5redmom}, to prove the upper bound in \eqref{r5hilnier}, it is enough to show that
\begin{multline}
\Ex\sup_{t\in B^m_2, x\in \kxx} \sum_{ijk} a_{ijk} x_i Y_j t_k \lesssim \ \Ex \sup_{t \in B^m_2} \sum_{ijk} a_{ijk} X_i Y_j t_k+\sup_{x \in \kxx} \Ex \sup_{t \in B^m_2}\sum_{ijk} a_{ijk} x_i Y_j t_k\\
+\sup_{t \in B^m_2} \norx{\left(\sqrt{\sum_j \left(\sum_k a_{ijk} t_k \right)^2}\right)_i}+\sup_{t\in B^m_2} \norxy{\left(\sum_k a_{ijk} t_k\right)_{ij}}. \label{r5p14}
\end{multline}

By \eqref{r5zawieranie} we have $\kxx\subset p^{1/2}B^n_2+pB^n_1$. We decompose
\[(p^{-1/2}\kxx)\times B^m_2=\bigcup_{l=1}^N (T_l+(z^l,s^l)),\ N\leq \exp(Cp)\] 
using Corollary \ref{r573}. Lemma \ref{r5510} yields 
\begin{multline}
\Ex\sup_{x \in \kxx, t\in B^m_2} \sum_{ijk} a_{ijk} x_i Y_j t_k\\
\lesssim  p^{1/2}\max_l \Ex\sup_{(x-z^l,t-s^l) \in T_l} \sum_{ijk} a_{ijk} x_i Y_j t_k+\sup_{t\in B^m_2} \norxy{(\sum_k a_{ijk}t_k)_{ij}} \label{r5p15}\\
\leq p^{1/2} \max_l\left(\Ex\sup_{(x,t) \in T_l} \sum_{ijk} a_{ijk} x_i Y_j t_k+\Ex\sup_{(x,t) \in T_l} \sum_{ijk} a_{ijk} x_i Y_j s^l_k+\Ex\sup_{(x,t) \in T_l} \sum_{ijk} a_{ijk} z^l_i Y_j t_k \right)\\
+\sup_{t\in B^m_2} \norxy{(\sum_k a_{ijk}t_k)_{ij}},
\end{multline}
where in the last line we used that $\Ex Y_j=0$.
Now because of the properties of our decomposition (obtained by Corollary \ref{r573}), and since $\Ex Y^4_j\leq C$  (recall \eqref{wnnorm}) 
\begin{equation}\label{r5p16}
\max_l \Ex\sup_{(x,t) \in T_l} \sum_{ijk} a_{ijk} x_i Y_j t_k\lesssim   p^{-1/2}\sqrt{\sum_{ijk} a_{ijk}^2}\lesssim p^{-1/2}\Ex \sup_{t \in B^m_2} \sum_{ijk} a_{ijk} X_i Y_j t_k, 
\end{equation}
where the last inequality follows by Lemma \ref{turbolemat}. By Corollary \ref{r573}   $(z^l,s^l) \in (p^{-1/2} \kxx )\times B^m_2$ so that $T_l \subset 2 (p^{-1/2} \kxx )\times 2 B^m_2$ and we may use Corollary \ref{r5some} to get that
\begin{multline}
\Ex\sup_{(x,t) \in T_l} \sum_{ijk} a_{ijk} x_i Y_j s^l_k\leq 2p^{-1/2} \sup_{t \in B^m_2} \Ex \sup_{x \in \kxx} \sum_{ijk} a_{ijk} x_i Y_j t_k\\
\lesssim p^{-1/2}\sup_{t\in B^m_2}\left( \norxy{\left(\sum_{k} a_{ijk}  t_k\right)_{ij} }+ \norx{ \left(\sqrt{\sum_j \left(\sum_k a_{ijk} t_k \right)^2} \right)_i } \right). \label{r5p18}
\end{multline}
Finally, since $z^l \in p^{-1/2}\kxx$,
\begin{equation}
\Ex\sup_{(x,t) \in T_l} \sum_{ijk} a_{ijk} z^l_i Y_j t_k\leq p^{-1/2}\sup_{x \in \kxx} \Ex \sup_{t \in B^m_2}\sum_{ijk} a_{ijk} x_i Y_j t_k. \label{r5p19}
\end{equation}
To prove \eqref{r5p14} it is enough to plug \eqref{r5p16}-\eqref{r5p19} into \eqref{r5p15}.
\end{proof}

\section{Expectation of suprema of a certain Gaussian process}\label{rgaus}

Let us fix a tensor $A=(a_{ijk})_{i,j\leq n,k\leq m}$. The main result of this section is Proposition \ref{proc}, in which we estimate the expectation of the supremum of a Gaussian process  $(G_{(x,t)})_{(x,t) \in V}$, where $V\subset p^{-1/2}\kxx \times T \subset (B_2+p^{1/2}B_1)\times \eR^m $, $G_{(x,t)}=\sum_{ijk} a_{ijk} g_i x_j t_k$ and $T$ is arbitrary. To estimate such a quantity, one needs to study the distance on $ p^{-1/2}\kxx \times T$ given by
$$d_A( (x,t),(x',t'))=(\Ex |G_{(x,t)}-G_{(x',t')}|^2)^{1/2}=\alpha_A(x\ten t - x' \ten t'),$$
where $x\ten t =(x_j \cdot t_k)_{j\leq n, k\leq m} \in \eR^{nm}$ and $\alpha_A$ is the norm (see Remark \ref{rem:norma} below) on $\eR^{nm}$ defined by the formula 
$$\alpha_A((x)_{jk})=\sqrt{\sum_i^n \left(\sum_{j=1}^n \sum_{k=1}^m a_{ijk} x_{jk} \right)^2}.$$
\begin{remark}\label{rem:norma}
We may assume that $\alpha_A$ is a norm (a priori it is only a seminorm).  We can replace the tensor $A$ by $\hat{A}=(\hat{a}_{ijk})_{ijk}$ where
\begin{equation}\label{eq:zmientensor}
\hat{a}_{ijk}=\begin{cases} a_{ijk} &i\leq n, \\ 0   &n<i\leq n+nm \textrm{ and } (j,k)\neq \sigma(i) \\ \eta &n<i\leq n+nm \textrm{ and } (j,k)=\sigma(i)  \end{cases},   
\end{equation}
$\sigma$ is any bijection between the sets $[nm]$ and $[n]\times [m]$, and $\eta>0$ is a small number.  Then $\alpha_{\hat{A}}$ is a norm given by
\[\alpha_{\hat{A}}((x))_{jk})=\sqrt{\sum_i^n \left(\sum_{j=1}^n \sum_{k=1}^m a_{ijk} x_{jk} \right)^2+\eta \lv x \rv_2}\stackrel{x\neq 0}{>}0.\]
It is enough to derive all the results for $\hat{A}$ and take $\eta\to 0$.   
\end{remark}

We use the scheme introduced in \cite{adlatmel}. To proceed, we need some entropy estimates for subsets of $(B_2+p^{1/2}B_1)\times \eR^m$. As usual $N(U,\rho,\eps)$ stands for the smallest number of closed balls with the diameter $\eps$ in metric $\rho$ that cover the set $U$.
The crucial idea is to consider the distribution of the vector $\eps(G_n+E_n)$, $\eps>0$  which is denoted by $\mu_{n,\eps}$ (we recall that $E_n$ is the symmetric exponential vector in $\eR^n$ independent of $G_n$).

Let $S\subset \eR^m$ and define the following norm on $\eR^n$ (a priori it is only a seminorm, but after applying Remark \ref{rem:norma} it is a norm)
\[\beta_{A,S}(x):=\Ex \sup_{s \in S} \left| \sum_{ijk} a_{ijk} g_i x_j s_k \right|.\]
By the classical Sudakov  minoration (Theorem \ref{sud1}), for every $x \in \eR^n$ there exists a set $S_{x,\eps} \subset S$  of cardinality at most $\exp(C\eps^{-2})$ such that 
\begin{align}\label{eq:locwar}
 \forall_{t\in S} \exists_{t' \in S_{x,\eps}} \ \alpha_A(x\te(t-t'))\leq \eps \Ex \sup_{s \in S} \left| \sum_{ijk} a_{ijk} g_i x_j s_k \right| = \eps\beta_{A,S}(x).   
\end{align}
We define the following measure on $ \eR^n \times S$:
\begin{equation}\label{eq:ref1}
  \hat{\mu}_{\eps,S}(C)=\int_{\eR^n} \sum_{t \in S_{x,\eps}}   \1_C((x,t))) d\mu_{n,\eps}(x).  
\end{equation}
The following technical lemma (together with standard measure theory considerations) ensures that there is no measurability problem under the integral in \eqref{eq:ref1}. 

\begin{lemma}\label{lem:mierz}
The sets $S_{x,\eps}\subset S$ (sets satisfying \eqref{eq:locwar}) can be defined in such way that $x\mapsto S_{x,\eps}$ is a simple function (i.e. this function has finitely many values and the preimage of each value is a Borel set).      
\end{lemma}
\begin{proof}
W.l.o.g. we can assume that for fixed, small $\eta>0$ we have  
\[
    \forall_{x\in \eR^n, s\in S\subset \eR^m} \sqrt{\sum_{i} \left(\sum_{jk} a_{ijk} x_j s_j \right)^2} \geq \eta \lv x \rv_2 \lv s \rv_2
\]
(see Remark \ref{rem:norma}). In particular for any vector $x\in \eR^n$
    \begin{align}\label{eq:locbret}
        \beta_{A,S}(x) &\geq \sup_{s\in S} \Ex |\sum_{ijk} a_{ijk}g_i x_j s_k|=\sup_{s\in S}  \sqrt{\sum_i \left(\sum_{jk} a_{ijk} x_j s_k \right)^2} \Ex |g_1| \geq \frac{\eta}{2}\lv x \rv_2 \sup_{s\in S} \lv s \rv_2.
    \end{align}
Since \eqref{eq:locwar} is homogeneous, it is sufficient to define $S_{x,\eps}$ on the unit sphere (and then put $S_{x,\eps}:=S_{x/\lv x \rv_2,\eps}$). We decompose the unit sphere $S_{n-1}$ into Borel sets such that
\[S_{n-1}=\bigcup_{U\in \mathcal{U}} U,\ \sup_{U\in \mathcal{U}} \Delta_2(U) \leq \frac{\eps \cdot \eta \sup_{s\in S} \lv s \rv_2}{(2+\eps)H(S)} \textrm{ and } |\mathcal{U}|<\infty,\ \]
  where $\Delta_2$ is the diameter of the set in the Euclidean metric, and
   \[H(S):=\sup_{t\in S}\sqrt{\sum_{ij}\left(\sum_k a_{ijk}s_k  \right)^2}+\Ex \sup_{s\in S} \sqrt{\sum_j \left(\sum_{ik} a_{ijk}g_i s_k \right)^2  }\in (0,\infty).\]
   Fix any $U\in \mathcal{U}$ and take any $x\in U$. Let $S_{x,\eps}\subset S$ be the set such that \eqref{eq:locwar} holds (for $x$). For any other $x'\in U$ we define $S_{x',\eps}:=S_{x,\eps}$. We will show that \eqref{eq:locwar} holds with $x$ replaced by $x'$ and with $2\eps$ instead of $\eps$. Pick any $t\in S$ and let $t'\in S_{x,\eps}=S_{x',\eps}$ be such that the following holds
   \begin{equation}\label{eq:locrow1}
    \alpha_A(x\te(t-t'))\leq \eps \Ex \sup_{s \in S} \left| \sum_{ijk} a_{ijk} g_i x_j s_k \right| = \eps\beta_{A,S}(x).   
   \end{equation}
   It suffices to show that (since $t\in S$ was arbitrary)
   \begin{align}\label{eq:locnaszcel}
    \alpha_A(x'\te(t-t'))\leq 2\eps\beta_{A,S}(x').   
   \end{align}
  Using the triangle inequality, the Cauchy-Schwarz inequality and \eqref{eq:locrow1} we obtain that
\[
   \alpha_{A}(x'\otimes (t-t'))\leq \alpha_{A}((x'-x)\otimes (t-t'))+\alpha_{A}(x\otimes (t-t')) \leq 2 H(S)\lv x-x'\rv_2+\eps\beta_{A,S}(x). 
\]
Again using the triangle inequality and the Cauchy-Schwarz inequality
\[\beta_{A,S}(x)\leq \beta_{A,S}(x')+\beta_{A,S}(x-x')\leq \beta_{A,S}(x')+\lv x - x' \rv_2 H(S).\]
The two inequalities above give
\begin{equation}\label{eq:loc21a}
 \alpha_{A}(x'\otimes (t-t'))\leq \eps \beta_{A,S}(x')+ \lv x-x'\rv_2 H(S)(2+\eps).   
\end{equation}
Using the upper bound on $\Delta_2(U)$ and \eqref{eq:locbret} (recall that $\lv x' \rv_2=1$) we get that
\begin{equation}\label{eq:loc21b}
 \lv x-x'\rv_2 H(S)(2+\eps)\leq \Delta_{2}(U)H(S)(2+\eps)\leq \eps \eta \sup_{s\in S} \lv s \rv_2 \leq \eps \beta_{A,S}(x').   
\end{equation}
Formulas \eqref{eq:loc21a} and \eqref{eq:loc21b} imply \eqref{eq:locnaszcel}.
\end{proof}

\begin{lemma}\label{entropia}
For any finite set $S\subset \eR^m$, $p>0$, $(x,t)\in (B^n_2+p^{1/2}B^n_1) \times S$ and $\eps>0$ we have
$$\hat{\mu}_{\eps,S}(B\bigl((x,t),d_A,r(\eps,x,t)\bigr))\geq \frac{1}{4}\exp(-\eps^{-2}/2-p^{1/2}\eps^{-1}), $$
where 
$$B\bigl((x,t),d_A,r(\eps,x,t)\bigr)=\{(x',t')\in \eR^n \times S: \alpha_A(x \ten t -x' \ten t') \leq r(\eps,x,t)\}$$ 
and
$$r(\eps,x,t)=C\left(\eps^2 \Ex\bb(E_n)+\eps \bb(x)+\eps \Ex \alpha_A(E_n \te t ) \right).$$
\end{lemma}

\begin{proof}
Let us fix $(x,t) \in B^n_2 \times S$ and $\eps>0$. Set
$$U=\left\{x' \in \eR^n:\ \bb(x')\leq C \eps \Ex \bb(E_n)+\bb(x),\ \alpha_A((x'-x)\te t)\leq C \eps \Ex \alpha_A(E_n \te t) \right\}.$$
For any $x'\in U$, there exists $t'\in S_{x',\eps}$ such that $\alpha_A(x' \te (t-t')) \leq \eps \bb(x')$.  By the triangle inequality
\[
\alpha_A(x\otimes t-x'\otimes t')\leq \alpha_A((x-x')\otimes t)+\alpha_A(x'\otimes(t-t'))\leq r(\eps,x,t).
\]
Thus, by Lemma \ref{pods}, 
\[\hat{\mu}_{\eps,S} \Bigl(B\bigl((x,t),d_A,r(\eps,x,t)\bigr)\Bigr) \geq \mu_{n,\eps}(U)\geq 1/4\exp(-\eps^{-2}/2-p^{1/2}\eps^{-1}).\]

\end{proof}

\begin{corollary}\label{rozbicie}
For any $p,\eps>0$, $V \subset (B^n_2+p^{1/2}B^n_1)\times S\subset \eR^n \times \eR^m$
\[
N\left(V,d_A,r(\eps)\right)
\leq 4\exp(C\ve^{-2}+Cp^{1/2}\eps^{-1}),\] 
where
\[r(\eps):=\ve^2\Ex\beta_{A,S}(E_n)
+\ve\sup_{(x,t)\in V}\beta_{A,S}(x)+\ve\sup_{(x,t)\in V}\Ex\alpha_A(E_n\otimes t)\approx \sup_{(x,t)\in V}r(\eps,x,t).\]
\end{corollary}

\begin{proof}
Let $N=N(V,d_A,r(\eps))$. Then there exist points $(x_i,t_i)_{i=1}^N$ in $V$ such that
\[d_A((x_i,t_i),(x_j,t_j))> r(\eps).\]
Note that the balls $B((x_i,t_i),d_A,r(\eps)/2)$ are disjoint and, by
Lemma \ref{entropia}, each of these balls has $\hat{\mu}_{\eps,S}$ measure at least $1/4 \exp(-C\eps^{-2}-Cp^{1/2}\eps^{-1})$. On the other hand we obviously have $\hat{\mu}_{\eps,S}(\eR^n\times S)\leq \exp(C\ve^{-2})$, which implies $N\leq 4\exp(C\ve^{-2}+Cp^{1/2}\eps^{-1})$.
\end{proof}

To make the notation more compact we define for $S\subset \eR^m$ and $V\subset \eR^n \times \eR^m$,
\begin{align*}
&s_A(S):=\Ex \sup_{s \in S} \left| \sum_{ijk} a_{ijk} g_{ij} s_k \right|+\Ex \sup_{s \in S} \left| \sum_{ijk} a_{ijk} g_i \wy_j  s_k \right|=\Ex \sup_{s \in S} \left| \sum_{ijk} a_{ijk} g_{ij} s_k \right|+\Ex\beta_{A,S}(E_n), \\
&F_A(V):= \Ex \sup_{(x,t) \in V} \sum_{ijk} a_{ijk} g_i x_j t_k,\\
&\da(V):=\mathrm{diam}(V,d_A)=\sup_{(x,t),(x',t') \in V} \alpha_A(x \ten t - x'\ten t').
\end{align*}

\begin{lemma}\label{rozbijT}
For any $S\subset \eR^m$ and $p\geq 1$ there exists a decomposition $S=\bigcup_{i=1}^N S_i$ such that $N\leq \exp(Cp)$ and for any $i\leq N$,
$$ \sup_{s,s' \in S_i}  \Ex\alpha_A(E_n\otimes (s-s'))\leq p^{-1/2} s_A(S).$$
\end{lemma}
\begin{proof}
It is enough to use the Sudakov minoration Theorem \ref{sud1} and observe that  Lemma \ref{turbolemat} implies that
\begin{equation}\label{rempor}
\Ex \alf(E_n \otimes t) \approx \sqrt{\sum_{ij} \left( \sum_k a_{ijk} t_k \right)^2}.    
\end{equation}
\end{proof}

\begin{lemma}\label{przesun}
Let $p\geq 1$, $V \subset(B^n_2+p^{1/2}B^n_1)\times S\subset \eR^n \times \eR^m$ and $(x,t)\in \eR^n \times \eR^m$. Then there exists a decomposition $V=\bigcup_{i=1}^N V_i$ such that $N\leq \exp(C2^{2l}p)$ and for any $i\leq N$
$$F_A(V_i+(x,t)) \leq F_A(V_i)+\bt(x)+C\Ex\alpha_A(E_n \ten t),$$
and
$$\da(V_i)\leq \frac{1}{2^{2l}p}s_A(S) + \frac{1}{2^l p^{1/2}} \sup_{(y,s) \in V}\left( \bt(y) +\Ex \alpha_A(E_n \ten s)  \right). $$
\end{lemma}
\begin{proof}
By Corollary \ref{r557} we decompose $(B^n_2+p^{1/2}B^n_1)=\bigcup_{i\leq N_1} U_i$, in such a way that $N_1\leq \exp(Cp)$ and
\begin{equation}
    \sup_{u,v\in U_i}\alpha_A((u-v) \ten t)\leq \frac{C}{p^{1/2}}\Ex \alpha_A(E_n \ten t).\label{loc1}
\end{equation}
Let $V_i=V\cap (U_i \times S)$.  Corollary \ref{rozbicie} with $\eps=2^{-l}p^{-1/2}$ yields the decomposition $V_i=\bigcup_{j\leq N_2} V_{ij}$, where $N_2\leq \exp(Cp2^{2l})$ and
\begin{align}
&\da(V_{ij}) \leq \frac{1}{2^{2l}p}s_A(S)+ \frac{1}{2^l p^{1/2}}\sup_{(y,s) \in V}  \left(\bt(y) + \Ex \alpha_A(E_n \otimes s) \right). \nonumber
\end{align}  

Since $\Ex \sum_{ijk} a_{ijk} g_i x_j t_k=0$, we have
$$F_A(V_{ij}+(x,t)) \leq F_A(V_{ij})+\bt(x)+\Ex \sup_{(y,s) \in V_{ij}} \sum_{i'j'k'} a_{i'j'k'} g_{i'} y_{j'} t_{k'}.$$
From Lemma \ref{r663}, and then by \eqref{rempor} and \eqref{loc1} we obtain
\begin{align*}
&\Ex \sup_{(y,s) \in V_{ij}} \sum_{i'j'k'} a_{i'j'k'} g_{i'} y_{j'} t_{k'} \\
&\lesssim \sqrt{\sum_{ij} \left( \sum_k a_{ijk} t_k \right)^2 }+p^{1/2} \sup_{u,v \in U_i}\sqrt{\sum_{i'}\left(\sum_{j'k'} a_{i'j'k'} (u_{j'}-v_{j'})t_{k'} \right)^2} \lesssim \Ex \alpha_A(E_n \ten t).
\end{align*}
Since $N_1N_2\leq \exp(C2^{2l}p)$, $V=\bigcup_{i\leq N_1,j\leq N_2} V_{ij}$ is the desired decomposition.
\end{proof}

\begin{corollary}\label{final}
Let $p\geq 1$, $T\subset \eR^m$ be arbitrary  and $V\subset (B^n_2+p^{1/2}B^n_1)\times S \subset \eR^n \times \eR^m $, where $S=T$ or $S=T-T$. Let also $V-V\subset (B^n_2+p^{1/2}B^n_1)\times (T-T)$ . Then there exists decomposition $V=\bigcup_{i=1}^N ((x_i,t_i)+V_i)$ such that $N\leq \exp(C2^{2l}p)$ and for each $i\leq N$
\begin{enumerate}[i)]
\item $(x_i,t_i)\in V$, $V_i-V_i \subset V-V,\ V_i \subset V-V $ and $\mathrm{card} (V_i)\leq \mathrm{card}(V)-1$,
\item $\sup_{(y,s)\in V_i} \left(\btt(y)+\Ex \alpha_A(E_n\ten s) \right) \leq  \frac{2}{2^lp^{1/2}}s_A(T)$, 
\item $\da(V_i) \leq \frac{1}{2^{2l}p}s_A(T)$,
\item $F_A(V_i+(x_i,t_i))\leq F_A(V_i)+2\btt(x_i)+C\Ex\alpha_A(E_n \ten t_i)$.
\end{enumerate}

\end{corollary}

\begin{remark}
The formulation of Corollary \ref{final} may seem unusual. The reason comes from our chaining argument. Our goal is to upper bound $F_A(V)$, where $V\subset (B_2+p^{1/2}B_1)\times T$. The first step is to decompose $V=\bigcup_i (V_i+(x_i,t_i))$ into smaller pieces using Corollary \ref{final}. But we cannot guarantee that $V_i\subset (B_2+p^{1/2}B_1)\times T$, only that $V_i\subset (V-V) \subset 2(B_2+p^{1/2}B_1)\times (T-T)$. Then we decompose each of the $V_i=\bigcup_j (x_{ij},t_{ij})+V_{ij}$ using this corollary, which ensures that $V_{ij}$ and $V_{ij}-V_{ij}$ are subsets of $V_i-V_i \subset 2(B_2+p^{1/2}B_1)\times (T-T)$. So after the first step, the boundary conditions stabilize. This is why one needs to consider both cases $S=T$ and $S=T-T$. It is also the reason why the definition of the coefficients $c(l,k)$ in the proof of Proposition \ref{proc} is different for $l=0$ and $l\geq 1$ (see below).
\end{remark}
\begin{proof}[Proof of Corollary \ref{final}]
By Lemma \ref{rozbijT} and Corollary \ref{r557} (applied for the norm $\bt(\cdot)$ with $\eps=2^{-l}p^{-1/2}$) we can find the decomposition's $S=\bigcup_{i=1}^{N_1} S_i, B^n_2+p^{1/2}B^n_1=\bigcup_{i=1}^{N_2} U_i$ such  that $N_1,N_2\leq \exp(C2^{2l}p)$ and
\begin{align*}
\sup_{s,s' \in S_i} \Ex \alpha_A(E_n \otimes (s-s')) &\leq \frac{1}{2^{l+1}p^{1/2}} \max\left( s_A(T),s_A(T-T)\right)\leq \frac{1}{2^lp^{1/2}}  s_A(T), \\
\sup_{x,x' \in U_i} \bt(x-x')&\leq \frac{1}{2^{l+1}p^{1/2}}s_A(S)\leq  \frac{1}{2^lp^{1/2}}  s_A(T) . 
\end{align*}
Let $V_{ij}=V\cap (U_i \times S_j)$. If $V_{ij} \neq \emptyset$ we take any point $(x_{ij},t_{ij}) \in V_{ij}$ and using Lemma \ref{przesun}  we decompose
$$V_{ij}-(x_{ij},t_{ij})=\bigcup_{k=1}^{N_3} V_{ijk}$$
in such a way that $N_3 \leq \exp(C2^{2l}p)$,
\begin{align*}
F_A(V_{ijk}+(x_{ij},t_{ij}))  &\leq F_A(V_{ijk})
+2\btt(x_{ij}) +C\Ex\alpha_A(E_n \otimes t_{ij})
\end{align*}
(trivially $\bt(\cdot) \leq 2\btt(\cdot)$) and
\begin{align*}
&\da(V_{ijk}) \leq \frac{1}{2^{2l}p} s_A(S) + \frac{1}{2^l p^{1/2}} \left(\sup_{(y,s) \in V_{ij}} \bt(y-x_{ij})+\sup_{(y,s) \in V_{ij}} \Ex \alpha_A(E_n \otimes (s-t_{ij})) \right)\\
&\leq \frac{1}{2^{2l}p}s_A(S) + \frac{1}{2^l p^{1/2}}\left(\sup_{y,y' \in U_i} \bt(y-y')+\sup_{s,s' \in S_j} \Ex \alpha_A(E_n \otimes (s-s')) \right)\lesssim  \frac{1}{2^{2l}p} s_A(T).
\end{align*}
(trivially $s_A(S)\leq 2 s_A(T)$). Observe that 
\begin{align*}
V_{ijk}-V_{ijk}\subset V_{ij}-(x_{ij},t_{ij})-(V_{ij}-(x_{ij},t_{ij}))\subset V_{ij}-V_{ij}\subset V-V
\end{align*}
and by an analogous argument, $V_{ijk}\subset V-V$.
The final decomposition is obtained by relabeling the decomposition $V=\bigcup_{ijk} ((x_{ij},t_{ij})+V_{ijk})$.

\end{proof}

\begin{proposition}\label{proc}
For any $p\geq 1$, any non-empty $T\subset \eR^m$ and  $V\subset (B^n_2+p^{1/2}B^n_1) \times T$, 
\begin{align*}
F_A(V) \lesssim  \frac{1}{p^{1/2}}s_A(T)+\sup_{(x,t) \in V} \btt(x)+\sup_{(x,t) \in V}  \Ex \alpha_A(E_n \otimes t)+p^{1/2}\da(V) .
\end{align*}
\end{proposition}
\begin{proof}
W.l.o.g we may assume that $V$ is finite and $V\subset (1/2 (B^n_2+p^{1/2}B^n_1)) \times T$, so that $V-V\subset (B^n_2+p^{1/2}B^n_1)\times (T-T)$. We define

\begin{align*}
\de_0:=\da(V),\quad \hd_0:=\sup_{(x,t) \in V}  \btt(x)+\sup_{(x,t) \in V} \Ex \alpha_A(E_n \otimes t),\\
\de_l:=2^{-2l}p^{-1}s_A(T), \quad     \hd_l:=2^{-l}p^{-1/2}s_A(T) \quad \textrm{for } l=1,2,\ldots.
\end{align*}
Let for $k=1,2,\ldots$
\begin{align*}
c(0,k):=\sup\{&F_A(U):  U\subset V,\ |U|\leq k\},\\
c(l,k):=\sup\{&F_A(U): U\subset V-V,\  U-U\subset V-V,\ |U|\leq k,\\
&\da(U)\leq \de_l,\ \sup_{(x,t)\in U} \left(\btt(x) + \Ex \alpha_A(E_n \otimes t) \right) \leq 2\hd_l \}\ \textrm{ for } l\geq 1.
\end{align*}
Clearly $c(l,1)=0$ and if $U\subset V$ then 
\[\da(U)\leq \de_0 \textrm{ and }  \sup_{(x,t)\in U} \left(\btt(x) + \Ex \alpha_A(E_n \otimes t) \right)\leq C\hd_0.\]
We will show that for $k>1$ and $l\geq0$ we have
\begin{equation}
c(l,k)\leq c(l+1,k-1)+C\left(2^{l}p^{1/2}\de_l+\hd_l \right). \label{cel}
\end{equation}
To this end take any set $U$ as in the definition of $c(l,k)$ (in particular $|U|=k$) and apply Corollary \ref{final} to it (with $l$ replaced by $l+1$) to obtain the decomposition $U=\bigcup_{i=1}^N ((x_i,t_i)+U_i)$. In particular (by Corollary \ref{final})
\begin{align*}
&U_i\subset V-V,\ U_i-U_i\subset V-V,\ \max_i |U_i| \leq |U|-1\leq k-1,\\
    &\max_{i\leq N} \da(U_i)\leq \frac{1}{2^{2(l+1)}p} s_A(T)=\de_{l+1},\\
  &\max_{i\leq N}  \sup_{(x,t)\in U_i} \left(\btt(x) + \Ex \alpha_A(E_n \otimes t) \right) \leq \frac{2}{2^{(l+1)}p^{1/2}} s_A(T)=2\hd_{l+1}. \ 
\end{align*}
Thus, the sets $(U_i))_{i\leq N}$ satisfy the conditions described in the definition of $c(l+1,k-1)$, so that
$$\max_{i \leq N} F(U_i) \leq c(l+1,k-1).$$
Lemma \ref{r5gauskonc}  yields
\begin{align*}
F_A(U)=F_A\Big( \bigcup_{i\leq N} U_i+(x_i,t_i)\Big)  \leq C \sqrt{\ln N} \da(U)+\max_i F(U_i+(x_i,t_i)).
\end{align*}
 Since $N\leq \exp(C2^{2l}p)$ (cf. Corollary \ref{final}) from the definition of $c(l,k)$, we obtain
$$\sqrt{\ln N} \da(U)\leq C 2^{l} p^{1/2} \de_l$$
and for each $i$ we have by Corollary \ref{final} (recall that $(x_i,t_i)\in U$)
\begin{align*}
F(U_i+(x_i,t_i))\leq F(U_i) +C\left(\btt(x_i)+\Ex \alpha_A(E_n \otimes t_i)\right)
\leq F_A(U_i)+C\hd_l.
\end{align*}
So we have proved \eqref{cel}. It implies that for any $k$ we have
\begin{align*}
c(l,k)&\lesssim \sum_{j=0}^\infty \left(2^{j/2}p^{1/2}\de_j+C\hd_j\right) \\
&\lesssim \frac{1}{p^{1/2}}s_A(T)+\sup_{(x,t) \in V}  \btt(x)+\sup_{(x,t) \in V}  \Ex \alpha_A(E_n \otimes t)+p^{1/2}\da(V). 
\end{align*}
As a result
$$F(V)\leq \sup_k c(0,k) \lesssim\frac{1}{p^{1/2}}s_A(T)+\sup_{(x,t) \in V}  \btt(x)+\sup_{(x,t) \in V}  \Ex \alpha_A(E_n \otimes t)+p^{1/2}\da(V) .$$
\end{proof}

\section{Expectation of suprema of a certain Exponential process }\label{rexp}
In this section, we derive an analog of Proposition \ref{proc} for exponential processes when $T=\Te$. We do not know how to derive such a result for the general set $T\subset \eR^m$. We follow the notation of Section \ref{rgaus}. We will often use the following identity which is valid thanks to Lemma \ref{turbolemat} 
\begin{align}\label{porsa}
s_A(\Te)\approx^q \sqrt[q]{\sum_k \left( \sum_{ij} a_{ijk}^2 \right)^{q/2}}.
\end{align}
We also introduce a new norm on $\eR^n$ (using Remark \ref{rem:norma} we can ensure that it is a norm), which is essential for the proof
\begin{equation}\label{defalfa1}
 \varphi_A(x):=\sqrt[2q]{\sum_k \left(\sum_i \frac{\left(\sum_j a_{ijk}x_j \right)^4}{\sum_j a_{ijk}^2} \right)^{q/2}}.   
\end{equation}
We will need a technical estimation involving the norm $\varphi_A$.
\begin{lemma}\label{osznorm}
Let $\varphi_A$ be the norm defined by \eqref{defalfa1}. Then
\[\Ex \varphi_A(E_n)\lesssim^q \sqrt{\wspol}.\]
\end{lemma}
\begin{proof}
By Theorem \ref{r5gl} (precisely formula \eqref{eq:refloc2}) we obtain
\[\Ex \left(\sum_j a_{ijk} \wy_j\right)^4\approx \left(\Ex \left(\sum_j a_{ijk} \wy_j\right)^2\right)^2=\left(\sum_j a_{ijk}^2\right)^2.\]
Fact \ref{hiper1} (applied to $\eR^n$ with the norm $\lv x \rv_4=\sqrt[4]{\sum_i x^4_i}$ and $r=2q$, $p=1$), Jensen's inequality, and the above give
\begin{align*}
 \Ex \left( \sum_i \frac{(\sum_j a_{ijk}\wy_j)^4}{\sum_j a^2_{ijk}}\right)^{q/2}&\lesssim^q \left(   \Ex \left(  \sum_i \frac{(\sum_j a_{ijk}\wy_j)^4}{\sum_j a^2_{ijk}} \right)^{1/4} \right)^{2q} \leq  \left( \sum_i    \Ex  \frac{(\sum_j a_{ijk}\wy_j)^4}{\sum_j a^2_{ijk}} \right)^{q/2} \\
 &\lesssim  \left(\sum_{ij} a^2_{ijk}\right)^{q/2}.   
\end{align*}
Again Jensen's inequality and the above  give
\begin{align*}
\Ex \varphi_A(E_n) \leq \sqrt[2q]{\sum_k \Ex \left( \sum_i \frac{(\sum_j a_{ijk}\wy_j)^4}{\sum_j a^2_{ijk}}\right)^{q/2}}\lesssim^q \sqrt[2q]{\sum_k \left(\sum_{ij} a^2_{ijk}\right)^{q/2}}\lesssim^q \sqrt{\wspol},
\end{align*}
where we used \eqref{porsa} in the last inequality .
\end{proof}

 The next lemma was inspired by \cite{adlat} (Theorem 7.2 therein).

\begin{lemma}\label{iteracja}
For any set $U\subset \eR^n$ we have that
\begin{align*}
\Ex^g \sup_{x \in U} \Ex^{g'} \sup_{t \in B_{q'}} \sum_{ijk} a_{ijk} g_i g'_i x_j t_k \lesssim^q \sqrt{\wspol}\sup_{x \in U}  \varphi_A(x). 
\end{align*}
\end{lemma}
\begin{proof}
Duality and Lemma \ref{turbolemat}  imply that
\begin{align}
\Ex^g \sup_{x \in U} \Ex^{g'} \sup_{t \in B_{q'}} \sum_{ijk} a_{ijk} g_i g'_i x_j t_k &\lesssim^q \Ex \sup_{x \in U} \sqrt[q]{\sum_k \left(\sum_i g^2_i \left(\sum_j a_{ijk} x_j \right)^2 \right)^{q/2}}.
\end{align}
By Cauchy-Schwarz's inequality applied  the summation over index $i$ 
\begin{align*}
&\Ex \sup_{x \in U} \sqrt[q]{\sum_k \left(\sum_i g^2_i \left(\sum_j a_{ijk} x_j \right)^2 \right)^{q/2}}\\
 &= \Ex \sup_{x \in U} \sqrt[q]{\sum_k \left(\sum_i g^2_i \sqrt{\sum_j a^2_{ijk}}\cdot  \frac{\left(\sum_j a_{ijk} x_j \right)^2}{\sqrt{\sum_j a^2_{ijk}} } \right)^{q/2}} \\
&\leq \Ex \sup_{x \in U} \sqrt[q]{\sum_k \left( \sqrt{\sum_{ij} g^4_i  a^2_{ijk}} \sqrt{\sum_i \frac{\left(\sum_j a_{ijk} x_j \right)^4}{\sum_j a^2_{ijk} } }\right)^{q/2}}=:H(U,A,q)\\
\end{align*}
Now by Cauchy-Schwarz's inequality applied  the summation over the index $k$ 
\begin{align*}
&H(U,A,q)\leq \Ex \sup_{x\in U}\sqrt[q]{\sqrt{\sum_k \left(\sum_{ij} g_i^4 a_{ijk}^2 \right)^{q/2}}\cdot \sqrt{\sum_k \left(\sum_i \frac{\left(\sum_j a_{ijk} x_j \right)^4}{\sum_j a_{ijk}^2} \right)^{q/2}}}\\
&=\sup_{x\in U} \varphi_A(x)\Ex\sqrt[2q]{\sum_k \left(\sum_{ij} g_i^4 a_{ijk}^2 \right)^{q/2}}   \lesssim^q  \sup_{x\in U} \varphi_A(x) \sqrt{\wspol},
\end{align*}
where the last line follows by Lemma \ref{obliczsq} and \eqref{porsa}.
\end{proof}

To make the presentation more compact, we introduce some new notation in the spirit of the previous section. For $y\in \eR^{n}\times \eR^m$, $x,x'\in \eR^n$, $t,t'\in \eR^m$ and $V\subset \eR^n\times \eR^m$ we denote
\begin{align*}
&\alpha_{\infty, A}((y_{jk})_{j\leq n, k\leq m}):=\max_i \left| \sum_{jk} a_{ijk} y_{jk} \right|, & &d_{\infty,A}((x,t),(x',t'))=\alpha_{\infty,A}(x\otimes t - x' \otimes t'),\\
&\Delta_{\infty,A}(V)=\mathrm{diam}(V,d_{\infty,A}).
\end{align*}
By using  Remark \ref{rem:norma} we can assume that all $\alpha_{\infty, A}$ is a norm.

\begin{fact}\label{exp}
Let $p\geq 1$, $V\subset (B_2\times+ p^{1/2}B_1)\times \Te$, and $\varphi_A(x)$  be the norm defined in \eqref{defalfa1}. Then 
\begin{multline}\label{exploc}
\Ex \sup_{(x,t)\in V} \sum_{ijk}a_{ijk} \wy_i x_j t_k \lesssim^q p^{-1/2}s_A(\Te)+\sup_{(x,t)\in V} \Ex \alf(E_n\otimes t)+\sqrt{\wspol}\sup_{x\in U}\varphi_A(x)\\
 +p^{1/2} \sup_{y\in B_2,(x,t)\in V}\alinf(y\otimes t) +p^{1/2} \Delta_{A}(V)
+p \diaminf(V).
\end{multline}
\end{fact}

\begin{proof}
 Lemma \ref{gauswyk}  yields that it is enough to show \eqref{exploc} with  
$\Ex \sup_{(x,t)\in V} \sum_{ijk}a_{ijk} \wy_i x_j t_k$ replaced by $ \Ex \sup_{(x,t)\in V} \sum_{ijk}a_{ijk} g_ig'_i x_j t_k$. Consider random tensor $A_g=(a_{ijk}g_i)_{ijk}$.
By applying  Proposition \ref{proc} conditionally on $g_1,g_2\ldots$  we get 
\begin{multline}\label{tuwstaw}
\Ex \sup_{(x,t)\in V} \sum_{ijk}a_{ijk} g_ig'_i x_j t_k \\
\lesssim \Ex^g \left(p^{-1/2}s_{A_g}(\Te)+\sup_{(x,t)\in V} \beta_{A_g,\Te}(x)+\sup_{(x,t)\in V}\Ex^{E_n} \alpha_{A_g}(E_n\otimes t)+p^{1/2}\Delta_{A_g}(V)\right),
\end{multline}
where $\Ex^{E_n}$ means expectation with respect to $E_n$ (conditionally on $g_1,\ldots,g_n)$.\\
By \eqref{porsa} and Lemma \ref{normlem}
\begin{equation}\label{a0}
\Ex s_{A_g}(\Te)\lesssim^q \Ex \sqrt[q]{\sum_k \left(\sum_{ij} g^2_i a^2_{ijk} \right)^{q/2}} \lesssim^q \wspol.
\end{equation}
Lemma \ref{iteracja} states that
\begin{equation}\label{apol}
\Ex^g\sup_{(x,t)\in V} \beta_{A_g,\Te}(x)=\Ex^{g} \sup_{(x,t) \in V} \Ex^{g'} \sup_{s \in \Te} \sum_{ijk} a_{ijk} g_i g'_i x_j s_k\lesssim^q \sqrt{\wspol}\sup_{(x,t) \in V} \varphi_A(x).
\end{equation}
The inequality \eqref{rempor}, Jensen's inequality  and Lemma \ref{por} give
\begin{align*}
&\Ex^g\sup_{(x,t)\in V}\Ex^{E_n} \alpha_{A_g}(E_n\otimes t)\lesssim \sqrt{\Ex \sup_{(x,t) \in V}  \sum_{ij} g_i^2\left( \sum_k a_{ijk} t_k \right)^2}\\
&\leq \sup_{(x,t) \in V}  \sqrt{\sum_{ij} \left( \sum_k a_{ijk} t_k \right)^2 }+\sqrt{\Ex \sup_{(x,t) \in V} \sum_{ij} (g_i^2-1)\left( \sum_k a_{ijk} t_k \right)^2}  \\
&\lesssim \sup_{(x,t) \in V}  \Ex  \alf(E_n\otimes t) +\sqrt{\Ex \sup_{(x,t) \in V} \sum_{ij} g_i g'_i\left( \sum_k a_{ijk} t_k \right)^2}.
\end{align*}
Since 
\[\sup_{y\in B_2,(x,t)\in V} \alinf (y\otimes t)=\max_i \sup_{t\in T} \sqrt{\sum_j \left( \sum_k a_{ijk} t_k \right)^2}, \]
  Lemmas \ref{kontrak}, \ref{turbolemat} and the AM-GM inequality imply 
\begin{align}
\sqrt{\Ex \sup_{(x,t) \in V} \sum_{ij} g_i g'_i\left( \sum_k a_{ijk} t_k \right)^2} \lesssim p^{1/2}  \sup_{y\in B_2,(x,t)\in \Te} \alinf (y\otimes t)+p^{-1/2}\wspol.\label{a2}
\end{align}
By Jensen's and triangle inequalities and Lemma \ref{por}
\begin{align*}
&\Ex \sup_{(x,t),(x',t')\in V}\sqrt{\sum_i g^2_i \left(\sum_{jk} a_{ijk} (x_j t_k - x'_j t'_k) \right)^2}\\
&\leq\sqrt{\left(\Delta_{A}(V)\right)^2+ \Ex \sup_{(x,t),(x',t')\in V}\sum_i (g^2_i-1) \left(\sum_{jk} a_{ijk} (x_j t_k - x'_j t'_k) \right)^2}\\
&\leq \Delta_{A}(V)+C\sqrt{\Ex \sup_{(x,t),(x',t')\in V} \sum_i g_i g'_i \left(\sum_{jk} a_{ijk} (x_j t_k - x'_j t'_k) \right)^2}.
\end{align*}
The above and Lemma \ref{uff} ($\diaminf(V) =\max_i \sup_{(x,t),(x',t')\in V} |\sum_{jk} a_{ijk} (x_j t_k - x'_j t'_k)| $) imply
\begin{align}
p^{1/2}\Ex \Delta_{A_g}(V)&=p^{1/2} \Ex \sup_{(x,t),(x',t')\in V}\sqrt{\sum_i g^2_i \left(\sum_{jk} a_{ijk} (x_j t_k - x'_j t'_k) \right)^2} \nonumber\\
&\leq p^{1/2}\Delta_{A}(V)+Cp^{1/2}\sqrt{\diaminf(V) \Ex \sup_{(x,t),(x',t') \in V}  \sum_i g_i g'_i  \sum_{jk} a_{ijk} (x_j t_k-x'_jt'_k) } \nonumber\\
&\leq   p^{1/2}\Delta_{A}(V)+Cp\diaminf(V) 
 +\frac{1}{4}\Ex \sup_{(x,t),(x',t') \in V}  \sum_i g_i g'_i  \sum_{jk} a_{ijk} (x_j t_k-x'_jt'_k)\nonumber\\
&\leq   p^{1/2}\Delta_{A}(V)+Cp\diaminf(V) +\frac{1}{2} \Ex \sup_{(x,t)\in V}  \sum_i g_i g'_i  \sum_{jk} a_{ijk} x_j t_k.\label{a3}
\end{align}
Applying \eqref{a0}-\eqref{a3} into \eqref{tuwstaw} yields

\begin{align*}
&\Ex \sup_{(x,t)\in V} \sum_{ijk}a_{ijk} g_ig'_i x_j t_k\\
&\leq C\Bigg(p^{-1/2}s_A(\Te)+\sup_{(x,t) \in V}  \Ex \alf(E_n\otimes t) +\sqrt{\wspol} \sup_{(x,t) \in V} \varphi_A(x)\\
&+ p^{1/2}  \sup_{y\in B_2,(x,t) \in V} \alinf(y\otimes t) +p^{1/2} \Delta_{A}(V)
+p \diaminf(V) \Bigg)+\frac{1}{2}\Ex \sup_{(x,t) \in V}  \sum_{ijk}   a_{ijk} g_i g'_i   x_j t_k.
\end{align*}
The assertion follows.
\end{proof}

\section{Main Decomposition Lemma}\label{rrozbij}
The goal of this section is to provide decomposition $(p^{-1/2}B^X_p)\times B_{q'}=\bigcup_{l=1}^N (x_l,t_l)+U_l$ such that $N$ is "not big" and $\max_{l\leq N}\Ex \sup_{(x,t)\in U_l} \sum_{ijk} a_{ijk}\wy_i x_j t_k$ is "small". Our approach is to break down the set $(p^{-1/2}B^X_p)\times B_{q'}$ in such a way that each term on the right-hand side in \eqref{exploc} is "small", except $p^{1/2}\sup_{x\in B_2,t\in T} \alinf(x\otimes t)$. We do not know how to deal with this parameter which is the main obstacle to a more general result.

Before we prove a new entropy estimate we state some additional facts. The idea of the first one is taken from \cite{some}.
\begin{lemma}\label{staref1}
Let $\lv \cdot \rv$ be any norm on $\eR^n$, $t\in \eR^n$. Consider the set
\begin{equation}\label{aeps}
S_\eps:=\{x\in \eR^n: \alinf(x\otimes t) \leq \eps^2 \Ex \alpha_A(G_n\otimes t),\ \lv x \rv \leq C\eps\Ex\lv E_n\rv \}.
\end{equation}
Then for any $y\in B_2$ (we recall that $\gamma_{n,\eps}$ is the distribution of the $\eps G=\eps(g_1,\ldots,g_n)$)
\[\gamma_{n,\eps}(S_\eps+y)\geq ce^{-\frac{c}{\eps^2}}.\]
\end{lemma}
\begin{proof}
By Lemma \ref{kontr}, $\Ex \lv G_n \rv \leq 3\Ex \lv E_n \rv$ so Chebyshev's inequality gives
\begin{equation}\label{elo1}
\gamma_{n,\eps}\left(\lv x \rv\leq C\eps \Ex \lv E_n\rv \right)\geq 1-\frac{3}{C}.
\end{equation}
By Royen's Theorem \ref{roy}
\begin{align*}
\gamma_{n,\eps}(S_\eps) &\geq \gamma_{n,\eps} \left(x\in \eR^n :\alinf(x\otimes t)\leq \eps^2 \Ex \alf(E_n\otimes t)\right) \gamma_{n,\eps}\left(x\in \eR^n: \lv x\rv \leq C\eps\Ex\lv E_n \rv \right)\\ &\geq
 e^{-\frac{2}{\eps^2}}(1-1/C)\geq c e^{-\frac{2}{\eps^2}},
\end{align*}
where the last  inequality follows from Lemma \ref{maxexp} (cf. \eqref{rempor}) and \eqref{elo1}.\\
Since $S_\eps$ is symmetric, we obtain for any $y\in B_2$
\begin{align*}
\gamma_{n,\eps}(S_\eps+y)&=e^{-\frac{|y|^2}{2\eps^2}}\int_{S_\eps} e^{\langle y,x \rangle /\eps^2}d\gamma_{n,\eps}(x)=e^{-\frac{|y|^2}{2\eps^2}}\int_{S_\eps} \frac{1}{2} \left( e^{\langle y,x \rangle / \eps^2}+ e^{-\langle y,x \rangle / \eps^2}  \right) d\gamma_{n,\eps}(x)\\
&\geq \exp(-\eps^{-2}/2)\gamma_{n,\eps}(S_\eps)\geq c \exp(-c \eps^{-2}).
\end{align*}
\end{proof}

Before stating the next fact we recall that $\nu_{n,\eps}$ stands for the distribution of $\eps E_n$.
\begin{lemma}\label{Wniosek2}
Let $\lv \cdot \rv$ be any norm on $\eR^n$, $t\in \eR^n$, $S_{\eps}$ set defined in \eqref{aeps}. Then for any $p>0$ and $y\in p^{1/2}B_1$ we have
\[\nu_{n,\eps}(S_\eps+y)\geq ce^{-\frac{c}{\eps^2}-\frac{p^{1/2}}{\eps}}.\]
\end{lemma}
\begin{proof}
Since $\Pro(E_n\in U\subset \eR^n)=\nu_{n,\eps}(U/\eps)$, Theorem \ref{roymiesz} implies
\begin{align*}
\nu_{n,\eps}(S_\eps) &\geq \nu_{n,\eps}\left(\alinf(x \otimes t)\leq \eps^2\Ex \alf(G_n\otimes t)\right)\nu_{n,\eps} (\lv x\rv \leq C\eps\Ex\lv E_n\rv )\\
&\geq \frac{1}{4}e^{\frac{-4}{\eps^2}}\left(1-\frac{1}{C}\right),
\end{align*}
where in the last inequality we used Lemma \ref{maxexp}  ($\Ex \alf(G_n\otimes t) \approx \sqrt{\sum_{ij} (\sum_k a_{ijk} t_k)^2} $ by \eqref{rempor}) and Chebyshev's inequality. So for any $y\in p^{1/2}B_1$ we have
\begin{align*}
\nu_{n,\eps}(S_\eps+y)&=(2\eps)^{-n}\int_{S_\eps} \exp\left(-\frac{1}{\eps} \sum_{i=1}^n |x_i+y_i| \right)dx\geq \exp\left(-\frac{1}{\eps} \sum_{i=1}^n |y_i| \right)\int_{S_\eps}d\nu_{n,\eps}\\
&\geq \exp(-p^{1/2}/\eps)\nu_{n,\eps}(S_\eps)\geq ce^{-\frac{c}{\eps^2}-\frac{p^{1/2}}{\eps}}.
\end{align*}
\end{proof}
Recall that $\mu_{n,\eps}$ is the distribution of $\eps(G_n+E_n)$, where $G_n,E_n$ are independent.

\begin{lemma}\label{Lemat5}
Let $\lv \cdot \rv$ be an arbitrary norm on $\eR^n$. Consider any $y\in B_2+p^{1/2}B_1$, where $p>0$. 
If $S_\eps$ is the set defined in \eqref{aeps} then 
\[\mu_{n,\eps}(S_\eps+y)\geq ce^{-\frac{c}{\eps^2}-c\frac{p^{1/2}}{\eps}}.\]
\end{lemma}
\begin{proof}
We have $y=y'+y''$, where $y'\in B_2$ and $y''\in p^{1/2}B_1$. It is easy to check that
\[\left(S_{\eps/2}+y'\right)+\left(S_{\eps/2}+y''\right)\subset S_\eps+y,\]
where the left sum is the Minkowski addition.
Since $\mu_{n,\eps}$ is a convolution of $\gamma_{n,\eps}$ and $\nu_{n,\eps}$, we obtain
\[\mu_{n,\eps}(S_\eps+y)\geq \mu_{n,\eps}(\left(S_{\eps/2}+y'\right)+\left(S_{\eps/2}+y''\right))\geq \gamma_{n,\eps}(S_{\eps/2}+y')\nu_{n,\eps}(S_{\eps/2}+y'').\]
It is enough to apply Lemmas \ref{staref1},\ref{Wniosek2} to the above.
\end{proof}
Now we will follow the reasoning from Section \ref{rgaus}, where we defined the pivotal measure $\hat{\mu}_{\eps,T}$.
Theorem \ref{latsudmin} implies that for any $x\in \eR^n$ there exists $T_{x,\eps}\subset \Te$ such that $|T_{x,\eps}|\leq \exp(C\eps^{-2})$ and

\begin{equation}\label{loc1111}
\forall_{t\in \Te} \exists_{t'\in T_{x,\eps}} \dinf((x,t),(x,t')) \leq\eps^2 \Ex \sup_{t\in \Te} \sum_{ijk} a_{ijk} \wy_i x_j t_k \approx^q \eps^2 \beet(x) .
\end{equation}
The last inequality follows from Lemma \ref{turbolemat}, which allows to replace the exponential variables by Gaussians (at the cost of a constant).

For $C\subset \eR^n \times \Te$ we define the measure
\[ \mi(C)=\int_{\eR^n}\sum_{t\in T_{x,\eps}} \1_C((x,t))d\mu_{n,\eps}(x). \]
One can prove a counterpart of Lemma \ref{lem:mierz} in the above case, so there are no problems with measurability in the definition of the measure $\mi$.
By the construction
\begin{equation}\label{ogrmiar}
 \mi(\eR^n \times \Te)\leq \exp(C\eps^{-2}).
\end{equation}

\begin{lemma}\label{Fakt6}
Let $p>0$ and $(x,t)\in (B_2+p^{1/2}B_1)\times \Te$. Consider
\[D_\eps((x,t))=\{(x',t')\in (B_2+p^{1/2}B_1)\times \Te:\ \dinf((x,t),(x',t')) \leq r(\eps,x,t)\},\]
where
\begin{equation}\label{defr}
r(\eps,x,t)=C(q)\left(\eps^3\wspol+\eps^2\beet(x)+C\eps^2\Ex \alf(E_n\otimes t)\right).
\end{equation}
Then
\[\mi(D_\eps((x,t)))\geq c e^{-\frac{c}{\eps^2}-\frac{cp^{1/2}}{\eps}}.\]
\end{lemma}
\begin{proof}
Define
\begin{align*}
U:=\{x'\in (B_2+p^{1/2}B_1):\dinf((x,t),(x',t)) \leq C\eps^2 \Ex \alf(E_n\otimes t),\\
\beet(x')\leq \eps \Ex \beet(E_n) +\beet(x)\}.
\end{align*}
Using \eqref{loc1111}
\begin{equation}\label{loca}
\forall_{x'\in U} \exists_{t'\in T_{x',\eps}}\dinf((x',t),(x',t')) \lesssim^q \eps^2 \beet(x').
\end{equation}
Fix $x'\in U$ and pick $t'$ as above. The triangle inequality implies
\begin{align*}
\dinf((x,t),(x',t')) &\leq \dinf((x',t),(x',t'))+\dinf((x',t),(x,t))\\
&\lesssim^q \eps^2 \beet(x')+\eps^2\Ex \alf(E_n\otimes t)\leq r(\eps,x,t).
\end{align*}
The second inequality follows by \eqref{loca}, the last by the definition of $U$ (Lemma \ref{turbolemat} implies that $\Ex \beet(E_n)\lesssim^q \wspol$). Thus, 
\[\mi\Big(D_\eps((x,t))\Big)\geq \mu_{n,\eps}(U) \geq c e^{-\frac{c}{\eps^2}-\frac{cp^{1/2}}{\eps}},\]
where in the last inequality we used Lemma \ref{Lemat5} (we pick the norm $\beet(\cdot)$ so that $S_\eps +x\subset U$).
\end{proof}

The crucial entropy bound is a standard consequence of Lemma \ref{Fakt6} (cf. proof of Corollary \ref{rozbicie} and recall \eqref{ogrmiar}). 
\begin{fact}\label{rozkladdmax}
Suppose $V \subset (B_2+p^{1/2}B_1)\times \Te \subset\eR^n \times \eR^m$, $p\geq 1$. Then
\[N(V, \dinf,\sup_{(x,t)\in V}r(p^{1/2},x,t))\leq \exp(Cp),\]
where $r(p^{1/2},x,t)$ is defined in \eqref{defr}.
\end{fact}

\begin{lemma}\label{lemdorozbijwykl}
Let  $p\geq 1$, $V\subset (B_2+p^{1/2}B_1)\times B_{q'}\subset \eR^n \times \eR^m$. Then there exists a decomposition $V=\bigcup_{i=1}^N V_i+(x_i,t_i)$ such that $N\leq e^{Cp}$ and for each $i\leq N$ we have
\begin{enumerate}
\item  $ (x_i,t_i)\in V$, 

\item\label{eloziom} $
\sup_{(x,t)\in V_i}\left(\Ex \alf(E_n\otimes t)+\beet(x)\right)\lesssim^q p^{-1/2}\wspol$,
\item \label{nafi} $\sup_{(x,t)\in V_i}\varphi_A(x)\lesssim^q p^{-1/2}\sqrt{\wspol}$, where $\varphi_A$ is the norm defined in \eqref{defalfa1},
 \item \label{eloyolo} $ \Delta_{A}(V_i) \lesssim^q p^{-1} \wspol$,
\item \label{yolo} $\diaminf(V_i)\lesssim^q p^{-3/2} \wspol.$
\end{enumerate}
\end{lemma}
\begin{proof}
 Corollary \ref{r557} and Lemma \ref{osznorm} imply that there exists a decomposition $B_2+p^{1/2}B_1=\bigcup_{i\leq N'} U_i$ such that $N'\leq e^{Cp}$, and
 \begin{equation}\label{dodd}
     \max_{i\leq N'} \sup_{x,y\in U_i} \varphi_A(x-y)\lesssim^q p^{-1/2}\sqrt{\wspol}.
 \end{equation}
Let $V_i=V\cap (U_i\times \Te)$. Take any $(x_i,t_i)\in V_i\subset V$ and let $V'_i=V_i-(x_i,t_i)$.  For each $i$ we take the decomposition $V'_i =\bigcup_{ij\leq e^{Cp}} (x_{ij},t_{ij})+V_{ij}$ obtained by Corollary \ref{final} with $l=1$. In particular $(x_{ij},t_{ij})\in V'_i$ and 
\begin{align}\label{rowloc}
\sup_{(x,t)\in V_{ij}}\left(\beet(x)+\Ex \alf(E_n\otimes t) \right)\leq Cp^{-1/2}\wspol,\ \Delta_A(V_{ij})\leq Cp^{-1}\wspol.
\end{align}
Using Fact \ref{rozkladdmax} we can now decompose each $V_{ij}=\bigcup_{k\leq e^{Cp}} V_{ijk}$ in such a way that
\begin{align*}
\max_{ijk} \diaminf(V_{ijk})&\leq Cp^{-3/2} \wspol+Cp^{-1} \sup_{(x,t)\in V_{ij}}\left(\beet(x)+\Ex \alf(E_n\otimes t) \right)\\
&\leq Cp^{-3/2} \wspol,
\end{align*}
where we used \eqref{rowloc} in the last inequality.
We show that $V =\bigcup_{ijk} (x_i,t_i)+(x_{ij},t_{ij})+V_{ijk}$ is the desired decomposition. Clearly \ref{yolo} holds. By Corollary \ref{final} and our construction $(x_{ij},t_{ij})\in V'_i=V_i-(x_i,t_i)$. Thus $(x_{ij},t_{ij})+(x_i,t_i)\in V_i\subset V$. Since $V_{ijk}\subset V_{ij}$, the inequality \eqref{rowloc} implies \ref{eloziom} and  \ref{eloyolo}.
Since $V_{ijk}\subset V'_i-(x_{ij},t_{ij})\subset V'_i-V'_i$ and $V'_i\subset V_i-V_i\subset (U_i-U_i)\times 2B_{q'}$, we have
\[\sup_{(x,t)\in V_{ijk}}\varphi_A(x)\leq 2\sup_{(x,t)\in V'_{i}}\varphi_A(x)\leq 2\sup_{x\in U_i-U_i}\varphi_A(x)\lesssim^q p^{-1/2}\sqrt{\wspol},\]
where we used \eqref{dodd}  in the last inequality. Thus, \ref{nafi} holds and the assertion follows.
\end{proof}

\begin{theorem}

\label{rozbijwyklmaly}
Let $p\geq 1$ and consider the set $B^X_p\subset \eR^n $ given by \eqref{kulax}. There exists a decomposition  $(p^{-1/2}B^X_p)\times B_{q'}=\bigcup_{l\leq N } (x_l,t_l)+V_l\subset \eR^n \times \eR^m$ such that $N\leq e^{Cp}$ and for each $l\leq N$
\begin{enumerate}
\item $(x_l,t_l)\in p^{-1/2}B^X_p\times B_{q'}$,
\item $ \Ex \sup_{(x,t)\in V_l} \sum_{ijk} a_{ijk} \wy_i x_j t_k \leq p^{-1/2} \wspol+p^{1/2}  \sup_{y\in B_2, (x,t)\in V_l}\alinf(y\otimes t).$
\end{enumerate}
\end{theorem}
\begin{proof}
By \eqref{r5zawieranie} $p^{-1/2}B^X_p\subset B_2+p^{1/2}B_1$ thus we may use Lemma  \ref{lemdorozbijwykl} to decompose $p^{-1/2}B^X_p\times B_{q'}=\bigcup_{l\leq e^{Cp} } ((x_l,t_l))+V_l$.
By Fact \ref{exp}
\begin{multline}\nonumber
\Ex \sup_{(x,t)\in V_l} \sum_{ijk}a_{ijk} \wy_i x_j t_k\lesssim^q p^{-1/2}s_A(\Te)+\sup_{(x,t)\in V_l} \Ex \alf(E_n\otimes t)+\sqrt{\wspol}\sup_{(x,t)\in V_l}\varphi_A(x)\\
 +p^{1/2} \sup_{y\in B_2,(x,t)\in V_l}\alinf(y\otimes t) +p^{1/2} \Delta_{A}(V_l)
+p \diaminf(V_l).
\end{multline}

The upper bound of $\max_l \Ex \sup_{(x,t)\in V_l} \sum_{ijk} a_{ijk} \wy_i x_j t_k$  is then obvious from the properties of our decomposition obtained from Lemma \ref{lemdorozbijwykl}.
\end{proof}

\section{Proofs of the main results}\label{r5s4}

We will begin by proving Theorem \ref{r5thmmain}. The reason why we need an additional assumption about subgaussanity is that our counterpart to Corollary \ref{r573}, Theorem \ref{rozbijwyklmaly}, involves an additional term that is generally suboptimal.  With this additional assumption, we can use Corollary \ref{final} instead.

\begin{proof}[Proof of Theorem \ref{r5thmmain}]
 Lemma \ref{r5redmom} reduces the statement of the theorem to the following inequality
\begin{align}\nonumber
&\Ex \sup_{x\in B^X_p} \lv \sum_{ij} a_{ij} x_i Y_j \rv_{L_q} \lesssim^q \gamma\Ex \lv \sum_{ij} a_{ij} X_i Y_j \rv_{L_q}+\sup_{x\in B^X_p} \Ex \lv \sum_{ij} a_{ij} x_i Y_j \rv_{L_q} \\
&+ \sup_{f \in \kul} \norx{\left(\sqrt{\sum_jf(a_{ij})^2}\right)_i}+\sup_{f \in \kul} \norxy{(f(a_{ij}))_{ij}}.\label{loc9}
\end{align}
Without loss of generality, we can assume that we sum up to $i,j\leq n$. We are only interested in the space $V=\mathrm{span} \left((a_{ij})_{i,j \leq n} \right)\subset L_q$. Since  $V$ is a finite-dimensional subspace of $L_q$, it embeds well in $\ell_q$. More precisely, by \cite[Chapter 19 Theorem 12]{Handbook}, there exists a linear operator $U:V\to \ell^m_q$, where $m$ is finite and $\lv U \rv \lv U^{-1} \rv \leq 2$, where $\lv \cdot \rv$ means the operator norm. So we may assume that $L_q=\ell^m_q$ with $m$ finite.  Then  $a_{ij}=(a_{ijk})_{k\leq m}$, $B^*(\ell^m_q)=B^m_{q'}$, and for any $x,y\in \eR^n$ we have
\[\lv \sum_{ij} a_{ij} x_i y_j \rv_{L_q}=\lv \sum_{ij} a_{ij} x_i y_j \rv_{\ell_q}=\sqrt[q]{\sum_k \left|\sum_{ij} a_{ijk} x_i y_j \right|^q}=\sup_{t\in \Te} \sum_{ijk} a_{ijk}x_iy_j t_k.\]
To prove \eqref{loc9} it is enough to show the following inequality
\begin{align}
&\Ex \sup_{x\in \kxx, t\in \Te} \sum_{ijk} a_{ijk} x_i Y_j t_k \lesssim^q \gamma \Ex \sup_{ t\in \Te} \sum_{ijk} a_{ijk} X_i Y_j t_k +\sup_{x\in \kxx} \Ex \sup_{ t\in \Te} \sum_{ijk} a_{ijk} x_i Y_j t_k  \nonumber \\
&+\sup_{t\in \Te} \norx{\left(\sqrt{\sum_j (\sum_k a_{ijk} t_k)^2}\right)_{i}}+\sup_{t\in \Te} \norxy{\left(\sum_k a_{ijk} t_k\right)_{ij}}. \label{loccel}
\end{align}
Because of \eqref{r5zawieranie} $B^X_p\subset p^{1/2}B_2+pB_1$, we can  decompose the set $p^{-1/2}B^X_p\times \Te$ using Corollary \ref{final} with $l=1$ so that $(p^{-1/2}\kxx )\times \Te =\bigcup_{l=1}^N ((x^l,t^l)+T_l)$, $N\leq \exp(Cp)$ and $T_l$ satisfies properties $i)-iv)$ from Corollary \ref{final}. By homogeneity and  Lemma \ref{r5510} we obtain
\begin{multline}
\Ex \sup_{x\in \kxx, t\in \Te} \sum_{ijk} a_{ijk} x_i Y_j t_k\\
=p^{1/2}\Ex \sup_{(z,s)\in (p^{-1/2}\kxx) \times \Te} \sum_{ijk} a_{ijk} z_i Y_j s_k \lesssim p^{1/2}\max_l \Ex \sup_{(z,s)\in (x^l,t^l)+T_l} \sum_{ijk} a_{ijk} z_i Y_j s_k \\
+\sup\left\{\sum_{ijk} a_{ijk} p^{1/2}z_i y_j s_k  :  s \in T,\sum_i \nx\left(p^{1/2}z_i\right)\leq p,\   \sum_j \my(y_j)\leq p \right\}\\
= p^{1/2}\max_l \Ex \sup_{(z,s)\in (x^l,t^l)+T_l} \sum_{ijk} a_{ijk} z_i Y_j s_k+ \sup_{s \in B_{q'}} \norxy{\left(\sum_k a_{ijk} s_k \right)_{ij}}. \label{r5loc44}
\end{multline}
 Since $\Ex Y_j=0$  for every $j$,  by the triangle inequality,
\begin{align}
\Ex \sup_{(z,s)\in (x^l,t^l)+T_l} \sum_{ijk} a_{ijk} z_i Y_j s_k&\leq \Ex \sup_{(z,s)\in T_l}  \sum_{ijk} a_{ijk} z_i Y_j s_k+ \Ex\sup_{(z,s)\in T_l}  \sum_{ijk} a_{ijk} x^l_i Y_j s_k \nonumber \\
&+ \Ex \sup_{(z,s)\in T_l}  \sum_{ijk} a_{ijk} z_i Y_j t^l_k. \label{r5loc0}
\end{align}
Corollary \ref{final} ensures that $(x^l,t^l) \in (p^{-1/2}\kxx )\times \Te$ and $T_l\subset  (p^{-1/2}\kxx-p^{-1/2}\kxx )\times (\Te-\Te) \subset 2\left(p^{-1/2}\kxx\times \Te \right)$  so
\begin{equation} \label{r5loc3}
\Ex\sup_{(z,s)\in T_l}  \sum_{ijk} a_{ijk} x^l_i Y_j s_k \lesssim p^{-1/2}\sup_{z \in \kxx} \Ex \sup_{s \in \Te} \sum_{ijk} a_{ijk} z_i Y_j s_k.
\end{equation} 
Since $T_l \subset 2(p^{-1/2}\kxx \times \Te )$, Corollary \ref{r5some}  yields

\begin{multline}
\Ex \sup_{(z,s)\in T_l}  \sum_{ijk} a_{ijk} z_i Y_j t^l_k \leq 2p^{-1/2}\sup_{s \in \Te} \Ex \sup_{z \in \kxx} \sum_{ijk} a_{ijk} z_i Y_j s_k \\
\lesssim p^{-1/2} \left(\sup_{s \in \Te} \norxy{\left( \sum_k a_{ijk} s_k \right)_{ij}}+\sup_{s \in \Te} \norx{\left( \sqrt{\sum_j \left(\sum_k a_{ijk} s_k \right)^2}\right)_i}\right). \label{r5loc4}
\end{multline}

Random variables $Y_1,\ldots, Y_n$ are subgaussian with constant $\gamma$ so Theorem \ref{r5sub} (it is easy to check that processes $( \sum_{ijk} a_{ijk} z_i Y_j s_k)_{(z,s) \in T_l}$ and $( \sum_{ijk} a_{ijk} z_i g_j s_k)_{(z,s) \in T_l}$ satisfy its assumptions)
\begin{equation}\label{porsubgaus}
\Ex \sup_{(z,s)\in T_l}  \sum_{ijk} a_{ijk} z_i Y_j s_k \lesssim   \gamma \Ex \sup_{(z,s)\in T_l}  \sum_{ijk} a_{ijk} z_i g_j s_k.
\end{equation}

Using Proposition \ref{proc} and the properties $i)-iv)$ of Corollary \ref{final} we obtain 
\begin{multline}
\Ex \sup_{(x,t)\in T_l} \sum_{ijk} a_{ijk} x_i g_j t_k \\
\lesssim p^{-1/2}s_A(\Te)+\sup_{(z,s) \in T_l} \beta_{A,\Te}(z)+\sup_{(z,s) \in T_l} \Ex \alpha_A(E_n \ten s)+p^{1/2}\da(T_l) \\
\lesssim p^{-1/2} \wspol \lesssim^q p^{-1/2} \Ex \sup_{t\in \Te} \sum_{ijk} a_{ijk} X_i Y_j t_k, \label{r5loc6}
\end{multline}
where we used Lemma \ref{turbolemat}  in the last inequality. The inequality \eqref{loccel} follows from \eqref{r5loc44}-\eqref{r5loc6}.
\end{proof}

\begin{proof}[Proof of Theorem \ref{Thmnewres}]
We proceed similarly to the proof of Theorem \ref{r5thmmain}. In this case, instead of \eqref{loccel} we have to prove that
\begin{align*}
\Ex \sup_{(x,t)\in B^X_p\times B_{q'}} \sum_{ijk} a_{ijk} x_i Y_j &t_k \lesssim^q \Ex \sup_{t\in \Te} \sum_{ijk} a_{ijk} X_i Y_j t_k+\sup_{x \in \kxx} \Ex  \su \sum_{ijk} a_{ijk} x_i Y_j t_k   \\
&+\sup_{t\in \Te} \norx{ \left(\sqrt{\sum_j (\sum_k a_{ijk}t_k)^2}\right)_i}+\sup_{t\in \Te} \norxy{(\sum_k a_{ijk}t_k)_{ij}} \\
&+p\max_i \sup_{t\in B_{q'}} \sqrt{\sum_j (\sum_k a_{ijk}t_k)^2}. 
\end{align*}
The first difference is that we decompose $p^{-1/2}\kxx \times \Te $ using Theorem \ref{rozbijwyklmaly}. The second (and last) difference is that \eqref{porsubgaus} is not necessarily true. But the assumptions of Theorem \ref{porexp} are satisfied so 
\begin{multline*}
\Ex \sup_{(z,s)\in T_l} \sum_{ijk} a_{ijk} z_i Y_j s_k \lesssim \Ex \sup_{(z,s)\in T_l} \sum_{ijk} a_{ijk} z_i \wy_j s_k \\
\lesssim^q p^{-1/2}\Ex \sup_{t\in \Te} \sum_{ijk} a_{ijk} X_i Y_j t_k +p^{1/2}\max_i \sup_{(x,t)\in T_l}\sqrt{\sum_j \left(\sum_k a_{ijk} t_k \right)^2},
\end{multline*}
where the last inequality follows from Theorem \ref{rozbijwyklmaly} (and we upper bound $\wspol$ using Lemma \ref{turbolemat}). The rest of the proof remains the same.
\end{proof}

\appendix

\section{Appendix}
In this section, we collect results from previous work used in this paper.

\begin{theorem}[Gluskin-Kwapie\'n estimate]\label{r5gl}
Let  $X_1,\ldots,X_n $ be independent, symmetric r.v.'s with LCT which fulfill normalization condition \eqref{normalizacja}. Then for any $p \geq 1, \ a_1,\ldots a_n \in \eR$  we have
\begin{align}\label{eq:refloc1}
 \lv \sum_i a_i X_i \rv_p \approx \norx{(a_i)_i}.
\end{align}
In particular for any $1\leq p,r <\infty$
\begin{equation}\label{eq:refloc2}
\lv \sum_i a_i X_i \rv_r \lesssim\max\left(1,\frac{r}{p}\right) \lv \sum_i a_i X_i \rv_p.    
\end{equation}
\end{theorem}
\begin{proof}
    The formula \eqref{eq:refloc1} was formulated slightly differently in \cite{glu}. The above formulation can be found in \cite{some} (Theorem 2 there). Since for $u\geq 1$ $\lv (a_i)_i \rv_{X,up}\leq u \lv (a_i)_i \rv_{X,p}$ (recall \eqref{r5skalowanie}), \eqref{eq:refloc2} is a consequence of \eqref{eq:refloc1}.
\end{proof}

\begin{lemma}\label{r5510}
Let  $X_1,\ldots,X_n $ be independent, symmetric r.v.'s with LCT satisfying the normalization condition \eqref{normalizacja}.
Then for any sets $T_1,\ldots,T_k \subset \eR^n$ and any $C\geq 1$
$$\Ex \sup_{t \in \bigcup_{l=1}^k T_l} \sum_i t_i X_i \lesssim \max_{l \leq k} \Ex \sup_{t \in T_l} \sum_i t_i X_i +  \sup_{t,t' \in \bigcup_{l=1}^k T_l }\lv{t-t'}\rv_{X,log(k)/C}.$$
\end{lemma}
\begin{proof}
For $C=1$ this was shown in \cite[Lemma $5.10$]{adlat}. The assertion follows by calling the formula \eqref{r5skalowanie}
\[\sup_{t,t' \in \bigcup_{l=1}^k T_l }\lv{t-t'}\rv_{X,log(k)}\leq C \sup_{t,t' \in \bigcup_{l=1}^k T_l }\lv{t-t'}\rv_{X,log(k)/C}.\]
\end{proof}

\begin{theorem}\label{r5lat1}
Let $a_1,\ldots,a_n \in F$ where $(F,\lv \cdot \rv)$ is a vector space with seminorm $\lv \cdot \rv$. Assume that $X_1,\ldots,X_n $ are independent, symmetric r.v.'s with LCT satisfying the normalization condition \eqref{normalizacja}. Let $B^*(F)$ be the unit ball in the dual space i.e.
\[B^*(F)=\left\{f \in F^\ast: \sup_{x\in F: \lv x \rv \leq 1} f(x) \leq 1 \right\}, \]
where $F^\ast$ is the linear space of all functionals on $F$ and define. Then for any $p \geq 1$   we have
\[ \lv \sum_i a_i X_i \rv_p \approx \Ex \lv \sum_i a_i X_i \rv + \sup_{f \in B^*(F)}\norx{(f(a_i))_i}=\Ex \lv \sum_i a_i X_i \rv + \sup_{x\in \kxx}\lv \sum_i a_i x_i \rv.\]
\end{theorem}
The theorem can be proved in the same way as was \cite[Theorem 1]{lat}. For the convenience of the reader, we provide a shorter argument.
\begin{proof}
By duality $\lv \sum_i a_i X_i \rv = \sup_{f \in B^*(F)}\sum_i f(a_i) X_i$. Thus the equality part in the theorem follows by interchanging the suprema
    \begin{align*}
      \sup_{f \in B^*(F)}\norx{(f(a_i))_i}=      \sup_{f \in B^*(F)} \sup_{x\in \kxx} \sum_i f(a_i) x_i   =\sup_{x\in \kxx}\lv \sum_i a_i x_i \rv.
    \end{align*}
We will now prove the "$\approx$" part. By duality and \cite[Theorem 2.3]{lattkocz} (recall \eqref{eq:refloc2})
\begin{align*}
 \lv \sum_i a_i X_i \rv_p &= \lv \sup_{f \in B^*(F)} \sum_i f(a_i) X_i \rv_p \lesssim \Ex \lv \sum_i a_i X_i \rv + \sup_{f \in B^*(F)} \lv  \sum_i f(a_i) X_i \rv_p\\
 &\lesssim \Ex \lv \sum_i a_i X_i \rv + \sup_{f \in B^*(F)}\norx{(f(a_i))_i},
\end{align*}
where we used \eqref{eq:refloc1}  in the last line . The reverse of the above inequality is obvious.
\end{proof}

\begin{fact}\label{hiper1}
Let  $X_1,\ldots,X_n $ be independent, symmetric r.v.'s with LCT. Let $a_i\in F$, where $(F,\lv \cdot \rv)$ is a vector space with a seminorm $\lv \cdot \rv$. Then for any $p,r\geq 1$
\[ \lv \sum_i a_i X_i \rv_r \leq \max\left(1, C\frac{r}{p} \right)\lv \sum_i a_i X_i \rv_p.\]
\end{fact}
\begin{proof}
By homogeneity, we may assume that $X_1,\ldots, X_n$ satisfy the normalization condition \eqref{normalizacja}. Then it is enough to apply Theorem \ref{r5lat1} and \eqref{eq:refloc2}.
\end{proof}

\begin{fact}\label{hiper}
Let $X_1,X_2,\ldots$, $Y_1,Y_2,\ldots$ be symmetric, independent r.v.'s with LCT. Let $a_{ij} \in F$, where $(F,\lv \cdot \rv)$ is a vector space with a seminorm $\lv \cdot \rv$. Then for any $p,r\geq 1$, 
$$\lv \sum_{ij} a_{ij} X_i Y_j \rv_{p} \approx^{r,p} \lv \sum_{ij} a_{ij} X_i Y_j \rv_r.$$ 
\end{fact}
\begin{proof}
\newcommand{\nor}[1]{\left| \! \left| \! \left| #1 \right| \! \right| \! \right|}
Applying Fact \ref{hiper1} conditionally, first time to a normed space $(F^\infty, \nor{\cdot})$, where $\nor{(a_i)_i}=\lv \sum a_i X_i \rv_{r}$ and then to $(F,\lv \cdot \rv)$  we conclude that
$$
\lv \sum_{ij} a_{ij} X_i Y_j \rv_{r} \approx^{r,p} \left(\Ex^Y  \left( \Ex^X \lv \sum_{ij} a_{ij} X_i Y_j \rv^{r} \right)^{p/r}  \right)^{1/p} \approx^{r,p} \lv \sum_{ij} a_{ij} X_i Y_j \rv_{p}.
$$

\end{proof}

\begin{lemma}\label{normlem}
Let  $X_1,X_2,\ldots$ be independent, symmetric r.v.'s with LCT  and satisfying the normalization condition \eqref{normalizacja}. Then for any real numbers $a_{ij}$ and $q\geq 1$ we have
\[\Ex \sqrt[q]{\sum_i \left(\sum_j a^2_{ij} X^2_j \right)^{q/2}}\approx^q \sqrt[q]{\sum_i \left( \sum_j a_{ij}^2 \right)^{q/2}}. \]
\end{lemma}
\begin{proof}
Fact \ref{hiper1} applied to a seminorm $\eR^n \ni x \mapsto \sqrt[q]{\sum_i \left(\sum_j a^2_{ij} x^2_j \right)^{q/2}} $ yields
\[\Ex \sqrt[q]{\sum_i \left(\sum_j a^2_{ij} X^2_j \right)^{q/2}}\approx^q  \sqrt[q]{\sum_i \Ex\left(\sum_j a^2_{ij} X^2_j \right)^{q/2}} \approx^q \sqrt[q]{\sum_i \left( \Ex \sum_j a^2_{ij} X^2_j \right)^{q/2}},  \]
where the latter "$\approx$" follows by applying Fact \ref{hiper1} for each $i\leq n$  to the norm $\eR^n \ni x\mapsto \sqrt{\sum_j a^2_{ij} x^2_j} $
 and $a_i=e_i$ ($e_1,\ldots,e_n$ is the standard orthonormal basis in $\eR^n$). Since \eqref{normalizacja}  implies $\Ex X^2_i \approx 1$ (recall \eqref{wnnorm}) we conclude the proof.
\end{proof}

\begin{lemma}\label{turbolemat}
Assume that $(X_{ij})_{ij},(Y_i)_i,(Z_j)_j$ are symmetric, independent r.v.'s with LCT. Assume also that they satisfy the normalization condition \eqref{normalizacja}. Then for any $q\geq 1$ and any real numbers $a_{ijk}$ we have
\begin{align} 
&\Ex \sqrt[q]{\sum_k \left|\sum_{ij} a_{ijk} X_{ij}Y_i Z_j \right|^q}\approx^q \sqrt[q]{\sum_k \left(\sum_{ij} a^2_{ijk} \right)^{q/2}},  \label{firstloc} \\
&\Ex \sqrt[q]{\sum_k \left|\sum_{ij} a_{ijk} Y_i Z_j \right|^q},\Ex \sqrt[q]{\sum_k \left|\sum_{ij} a_{ijk} X_{ij}Y_i \right|^q}\approx^q \sqrt[q]{\sum_k \left(\sum_{ij} a^2_{ijk} \right)^{q/2}}, \label{ndloc}\\
&\Ex \sqrt[q]{\sum_i \left|\sum_{j} a_{ij}X_j  \right|^q}\approx^q \sqrt[q]{\sum_i \left(\sum_{j} a^2_{ij} \right)^{q/2}}.\label{rdloc}
\end{align}
\end{lemma}
\begin{proof}
Since $(X_{ij})_{ij}$ satisfy \eqref{normalizacja}, by \eqref{wnnorm} $\Ex X^2_{ij}\approx 1$. By conditionally  using Fact \ref{hiper1} twice (first in $\ell_q$ then in $\eR$)  we get
\begin{align*}
&\Ex \sqrt[q]{\sum_k \left|\sum_{ij} a_{ijk} X_{ij}Y_i Z_j \right|^q}\approx^q \Ex^{Y,Z} \sqrt[q]{\Ex^X \sum_k \left| \sum_{ij} a_{ijk} X_{ij} Y_i Z_j\right|^q } \\
&\approx^q  \Ex^{Y,Z}\sqrt[q]{ \sum_k \left(\Ex \left|\sum_{ij} a_{ijk} X_{ij}Y_i Z_j \right|^2 \right)^{q/2}}\approx \Ex \sqrt[q]{ \sum_k \left(\sum_{ij} a^2_{ijk}Y^2_i Z^2_j  \right)^{q/2}} . 
\end{align*}
So \eqref{firstloc} follows by invoking Lemma \ref{normlem} twice (first conditionally on $Y$). \\
Formulas \eqref{ndloc},\eqref{rdloc} can be proved similarly as \eqref{firstloc}.
    \end{proof}

\begin{lemma}\label{obliczsq}
For any real numbers $a_{ijk}$ and any $q\geq 1$ we have
    $$\Ex\sqrt[2q]{\sum_k \left(\sum_{ij} a_{ijk}^2 g_i^4  \right)^{q/2}} \lesssim^q \sqrt[2q]{\sum_k \left(\sum_{ij}  a_{ijk}^2 \right)^{q/2}}.$$ 
\end{lemma}
\begin{proof}
Fix $k \in \eN$. By applying Fact \ref{hiper1} using seminorm $\lv x\rv=(\sum_{ij}  a^2_{ijk}x^4_i)^{1/4}$ we get
\[\Ex \left(\sum_{ij}a^2_{ijk} g^4_i  \right)^{q/2}=\Ex \left(\left(\sum_{ij} a^2_{ijk} g^4_i \right)^{1/4}\right)^{2q}\lesssim^q  \left(\Ex \left(\sum_{ij} a^2_{ijk} g^4_i \right)^{1/4}\right)^{2q}.   \]
So by Jensen's inequality, we get
\[ \Ex \left(\sum_{ij} a^2_{ijk} g^4_i \right)^{q/2} \lesssim^q \left( \left(\Ex \sum_{ij} a^2_{ijk} g^4_i \right)^{1/4}\right)^{2q}\lesssim  \left(\sum_{ij} a^2_{ijk} \right)^{q/2}. \]
The assertion follows since by Jensen's inequality
\[\Ex\sqrt[2q]{\sum_k \left(\sum_{ij}a_{ijk}^2 g_i^4  \right)^{q/2}} \leq \sqrt[2q]{\sum_k \Ex  \left(\sum_{ij}a_{ijk}^2  g_i^4 \right)^{q/2}}. \]
\end{proof}

\begin{theorem}\cite[Theorem 1]{some}\label{r5some1}
Let  $X_1,\ldots,X_n $, $Y_1,\ldots,Y_n$  be independent, symmetric r.v.'s with LCT satisfying the normalization condition \eqref{normalizacja}. Then for any real numbers $a_{ij}$ and any $p \geq 1$ we have
$$\lv \sum_{ij} a_{ij} X_i Y_j \rv_p \approx \norxy{(a_{ij})_{ij}}+\norx{\left(\sqrt{\sum_j a^2_{ij}} \right)_i} +\nory{\left(\sqrt{\sum_i a^2_{ij}} \right)_j}. $$
\end{theorem}

\begin{corollary}\cite[Corollary 3]{some}\label{r5some}
Let  $X_1,\ldots,X_n $, $Y_1,\ldots,Y_n$ be independent, symmetric r.v.'s with LCT satisfying the normalization condition \eqref{normalizacja}. Then for any real numbers $(a_{ij})_{ij}$ and any $p \geq 1$ we have
\begin{align*}
\Ex \nory{\left(\sum_i a_{ij} X_i \right)_j}\lesssim \left(\norxy{(a_{ij})_{ij}} + \nory{\left( \sqrt{\sum_i a^2_{ij} } \right)_j}  \right).
\end{align*}
\end{corollary}
\begin{theorem}\cite{royen, latroy}\label{roy}
Let $(G,G')$ have a joint  Gaussian and centered distribution in $\eR^{n+m}$. Then for any symmetric convex sets $K,L$ in $\eR^n, \eR^m$ respectively  we have
\[\Pro(G\in K, G' \in L) \geq \Pro(G\in K)\Pro(G'\in L).\] 
\end{theorem}

\begin{theorem}\label{roymiesz}
Let $\eps>0$ and $\nu_{n,\eps}$ be the distribution of $\eps E_n$. Then for any symmetric, convex sets $K,L$ in $\eR^n$  we have
\[\nu_{n,\eps}(K\cap L)\geq \nu_{n,\eps}(K)\nu_{n,\eps}(L).\]
\end{theorem}
\begin{proof}
We say that a r.v. $X$ is a Gaussian mixture (by definition) if there exist independent r.v.'s $Y,g$ such that $\Pro(Y<0)=0$, $g\sim \mathcal{N}(0,1)$ and $X,Yg$ have the same distribution. Let $X$ be symmetric, exponential r.v. ($X$ has p.d.f. $g(x)=1/2 \exp(-|x|)$). It can be checked that $\eps X$ has the same distribution as $\eps \sqrt{2|X|}g$ (see the remark $(i)$ after the proof of Lemma 23 in \cite[Theorem 16]{nayar}). Thus $\nu_{n,\eps}$ satisfies the assumptions of \cite[Theorem 17]{nayar}.
\end{proof}

\begin{lemma} \cite[Lemma 5.6]{adlat}\label{kontr}
For any norm $\lv \cdot \rv$ on $\eR^n$, $\Ex \lv G_n \rv \leq 3\Ex \lv E_n\rv$.
\end{lemma}

\begin{lemma}\label{pods}
Consider any $\eps,p>0$, any norms  $\alpha_1,\alpha_2$ on $\eR^n$ and $y\in B^n_2+p^{1/2}B^n_1$. Let
\[K=\{x\in \eR^n: \alpha_1(x-y)\leq C \eps \Ex \alpha_1(E_n),\ \alpha_2(x)\leq C\eps \Ex \alpha_2(E_n)+\alpha_2(y)\}.\]
Then we have
$$\mu_{n,\eps}(K)\geq 1/4 \exp\left(-\eps^{-2}-p^{1/2}\eps^{-1}\right), $$
where  $\mu_{n,\eps}$ is the distribution of $\eps(G_n+E_n)$.
\end{lemma}
\begin{proof}
  It is a consequence of \cite[Lemma 5.3]{adlat} with $s=t=\eps$ and Lemma \ref{kontr}.  
\end{proof}

\begin{lemma}\label{maxexp}
Let  $\hat{\mu}$ be  the distribution of $\eps G_n=\eps(g_1,\ldots,g_n)$ or $\eps E_n=\eps(\wy_1,\ldots,\wy_n)$. Then for any real-valued matrix $(a_{ij})$ we have
\[\hat{\mu}\left(x\in \eR^n: \max_j|\sum_{i=1}^n a_{ij}x_i |\leq \eps^2 \sqrt{\sum_{ij} a_{ij}^2}\right)\geq c \exp(-\frac{c}{\eps^2}). \]
\end{lemma}
\begin{proof}
 Clearly, $g_1,\ldots,g_n,\wy_1,\ldots,\wy_n$ are symmetric unimodal (a random variable is symmetric unimodal if it has a density with respect to the Lebesgue measure, which is symmetric and non-increasing on $[0,\infty)$ c.f. \cite{some}). 
So \cite[Lemma 4]{some} implies
\[\hat{\mu}\left(x\in \eR^n: \max_j|\sum_{i=1}^n a_{ij}x_i |\leq t\right)\geq \frac{1}{4} e^{-8\eps^2 \sum_{ij}\frac{a^2_{ij}}{t^2} }. \]
It is enough to take $t= \eps^2 \sqrt{\sum_{ij} a_{ij}^2}$.
\end{proof}

\begin{corollary}[Corollary $5.7$ from \cite{adlat} with $d=1$]\label{r557}
Let  $\alpha( \ \cdot \ )$ be a norm on $\eR^n$ and $\rho_\alpha$ be a distance on $\eR^n$ defined by $\rho_\alpha(x,y)=\alpha(x-y)$. Then for any $p>0, \eps \in (0,1]$,
$$N(B^n_2+p^{1/2}B^n_1,\rho_\alpha,C\eps \Ex \alpha(\mathcal{E}_1,\ldots,\mathcal{E}_n ))\leq \exp\left( \eps^{-2}+p^{1/2}\eps^{-1}  \right).$$ 
\end{corollary}

\begin{lemma}\cite[Lemma 6.3]{adlat}\label{r663}
For any real-valued matrix $(a_{ij})_{i,j\leq n}$, $p\geq 1$ and $U\subset B^n_2+p^{1/2}B^n_1$ we have
$$\Ex \sup_{x \in U} \sum_{ij} a_{ij} x_i g_j \lesssim \sqrt{\sum_{ij} a^2_{ij}}+p^{1/2} \cdot \sup_{x,x' \in U} \sqrt{\sum_j \left(\sum_i a_{ij} (x_i-x'_i) \right)^2} .$$
\end{lemma}

\begin{corollary}\cite[Corollary 7.3]{adlat}\label{r573}
Let $A=(a_{ijk})_{ijk}$ be a real-valued tensor and $Z_1,\ldots,Z_n$ be independent mean zero r.v.'s. Then for any $p \geq 1$ and $ T \subset (B^n_2+p^{1/2}B^n_1) \times (B^n_2+p^{1/2}B^n_1)  $ there exists a decomposition $T = \bigcup_{l=1}^N (x^l,y^l)+T_l $ such that $N\leq \exp(Cp)$, $(x^l,y^l) \in T$ and for every l,
$$\Ex \sup_{(x,y)\in T_l} \sum_{ijk} a_{ijk} Z_i x_j y_k \lesssim p^{-1/2} \sqrt{\sum_{ijk} a_{ijk}^2} \max_i \lv Z_i \rv_4.$$
\end{corollary}

\begin{lemma}\cite[Lemma 3]{latgaus}\label{r5gauskonc}
Let $(G_t)_{t \in T}$ be a centered Gaussian process and $T = \bigcup_{l=1}^m T_l$. Then
$$\Ex \sup_{t \in T} G_t \leq \max_{l \leq m} \Ex \sup_{t \in T_l} G_t + C\sqrt{\ln(m)} \sup_{t,t' \in T} \sqrt{\Ex (G_t-G_{t'})^2}.$$
\end{lemma}

\begin{theorem}\cite[Theorem 3.18 (Sudakov minoration)]{proinban}\label{sud1}
Let $T\subset \eR^n$ be arbitrary and  $d_2$ be the standard Euclidean distance. Then for any $\eps>0$
\[ N(T,d_2,\eps\Ex \sup_{t\in T} \sum_{i=1}^n t_i g_i ) \leq e^{\frac{C}{\eps^2}} .\] 
\end{theorem}

\begin{theorem}\label{latsudmin}
Let $T\subset \eR^n$ be arbitrary and consider $d_{\infty}(s,t):=\max_i |s_i-t_i|$. Then for any $\eps>0$ 
\[ N(T,d_{\infty},\eps \Ex \sup_{t\in T} \sum_{i=1}^n \exe_i t_i) \leq e^{\frac{C}{\eps}}.\]
\end{theorem}
\begin{proof}
Take $T'\subset T$, which is a maximal $\eps$ net with respect to the distance $d_{\infty}$, so that 
\begin{enumerate}
\item for any $s,t\in T'$ $d_{\infty}(s,t)> \eps$,
\item for any $t\in T$ there exists $t'\in T'$ such that $d_{\infty}(t,t')\leq \eps$.
\end{enumerate}
Let $N=|T'|$. By standard arguments
\begin{align}\label{entropialoc}
N(T,d_{\infty},2\eps)\leq N \leq N(T,d_{\infty},\eps).
\end{align}
Fix $s,t \in T'$. Jensen's inequality implies (since for any $i$, $\Ex \exe_i=0$)
\[ \lv \sum_i (t_i-s_i)\exe_i \rv_{\ln N} \geq \max_i \lv (t_i-s_i) \exe_i \rv_p = d_{\infty}(s,t) \left(\Gamma(1+\ln N)\right)^{\frac{1}{\ln N}} \geq C^{-1} \eps \ln N.\]
Using \cite[Theorem 1]{latsud} (cf. \cite[Theorem 1.3]{lattkocz})
\[2\Ex \sup_{t\in T} \sum_i t_i\exe_i\geq \Ex \sup_{s,t\in T} \sum_i (t_i-s_i)\exe_i \geq \Ex \sup_{s,t\in T'} \sum_i (t_i-s_i)\exe_i \geq C^{-1}\eps \ln N.\]
Using \eqref{entropialoc} and the above
\[N(T,d_{\infty},2\eps) \leq N\leq \exp\left(\frac{2C\Ex \sup_{t\in T}\sum_i t_i\exe_i}{\eps}\right).\]
It is enough to substitute $2\eps= \eps ' \Ex \sup_{t\in T}\sum_i t_i\exe_i$.
\end{proof}

\begin{theorem}\cite[Theorem 12.16]{proinban}\label{r5sub}
Let $(G_t)_{t \in T}$ be a Gaussian process and $(Y_t)_{t \in T}$ be a process such that for any $\lambda\in \eR$
\[\forall_{t,t'\in T}\  \Ex \exp(\lambda (Y_s-Y_t) )\leq \exp\left(\frac{\lambda^2}{2} \lv G_t - G_{t'} \rv_2^2\right).  \]
Then
$$\Ex \sup_{t \in T} Y_t \lesssim \Ex \sup_{t \in T} G_t.$$
\end{theorem}

\begin{lemma}\label{XXX}
Let $(X_t)_{t\in T}$ be a symmetric process. Then for any fixed $t_0\in T$
\[\Ex \sup_{t\in T} |X_t|\leq 2\Ex \sup_{t\in T} X_t+\Ex |X_{t_0}|.\]
\end{lemma}
\begin{proof}
Clearly
\[\Ex \sup_{t\in T} |X_t|\leq \Ex \sup_{t\in T} |X_t-X_{t_0}|+\Ex|X_{t_0}|\leq \Ex \sup_{s,t\in T} |X_t-X_s|+\Ex|X_{t_0}|.\]
By the symmetry
\[\Ex \sup_{s,t\in T} |X_t-X_s|=\Ex \sup_{s,t\in T} (X_t-X_s)=\Ex\sup_{t\in T} X_t +\Ex \sup_{s\in T}(-X_s)=2\Ex\sup_{t\in T} X_t.\]
\end{proof}

\begin{lemma}\label{konteps}
   Consider $T\subset \eR^n$ and $x,y\in \eR^n$ such that for any $i\leq n$, $|x_i|\leq |y_i|$. Then
   \[\Ex \sup_{t\in T} \sum_i t_i x_i \eps_i \leq \Ex \sup_{t\in T}\sum_i t_i y_i \eps_i. \]
\end{lemma}
\begin{proof}
Since $\varphi_i(t)=\1_{y_i\neq 0} t \frac{x_i}{y_i}$ is a contraction, the assertion follows from \cite[Theorem 6.5.1]{tal}.
\end{proof}

\begin{lemma}\label{por}
For any set $T\subset \eR^n$
\[\Ex \sup_{t\in T} \sum_i t_i (g^2-1)\leq 2\Ex \sup_{t\in T} \sum_i t_i g_ig'_i.\]
\end{lemma}
\begin{proof}
If $g,g'$ are independent $\mathcal{N}(0,1)$ r.v.'s then $g-g',g+g'$ are independent $\mathcal{N}(0,2)$ r.v.'s. So by Jensen's inequality
\begin{align*}
\Ex \sup_{t\in T} \sum_i t_i (g^2_i-1)&\leq \Ex \sup_{t\in T} \sum_i t_i (g^2_i-(g')^2_i)= 2\Ex \sup_{t\in T} \sum_i t_i g_ig'_i.
\end{align*}
    \end{proof}

\begin{fact}\label{porexp}
Assume that the symmetric r.v.'s $X_1,X_2,\ldots$ are independent, have LCT, and satisfy the normalization condition \eqref{normalizacja}. Then for any set $T\subset \eR^n$ we have
\[\Ex \sup_{t\in T} \sum_{i=1}^n t_i X_i \leq 2 \Ex \sup_{t\in T} \sum_{i=1}^n t_i\wy_i.\]
\end{fact}
\begin{proof}
Consider Rademacher r.v.'s $\eps_1,\eps_2,\ldots $  which are independent of $X_1,X_2,\ldots$. From \eqref{r5cos}, $\Pro(|X_i|\geq t) \leq \Pro(|\wy_i|\geq t)$ for $t\geq 1$. So we may assume that $|X_i|\1_{|X_i|\geq 1} \leq |\wy_j|$ (by inversing the CDF on the changed probability space). Thus, by the symmetry of $X_1,X_2,\ldots$ and Lemma \ref{konteps} (applied conditionally on $X_i,\wy_i$) we get
\[\Ex \sup_{t\in T} \sum_{i=1}^n t_i X_i\1_{|X_i|\geq 1}=\Ex \sup_{t\in T} \sum_{i=1}^n t_i |X_i|\eps_i \1_{|X_i|\geq 1} \leq \Ex \sup_{t\in T} \sum_{i=1}^n t_i |\wy_i|\eps_i =\Ex \sup_{t\in T} \sum_{i=1}^n t_i\wy_i.\]
Using a similar argument and Jensen's inequality
\begin{align*}
\Ex \sup_{t\in T} \sum_{i=1}^n t_i X_i\1_{|X_i|< 1}&=\Ex \sup_{t\in T} \sum_{i=1}^n t_i |X_i|\eps_i \1_{|X_i|< 1} \leq \Ex \sup_{t\in T} \sum_{i=1}^n t_i \eps_i \\
&=\Ex \sup_{t\in T} \sum_{i=1}^n t_i \eps_i(\Ex |\wy_i|)\leq \Ex \sup_{t\in T} \sum_{i=1}^n t_i \eps_i |\wy_i|=\Ex \sup_{t\in T} \sum_{i=1}^n t_i \wy_i.
\end{align*}

\end{proof}

\begin{lemma}\label{gauswyk}
Let $(g_i)_i,(g'_i)_i$ be independent $\mathcal{N}(0,1)$ r.v.'s. Then for any $T\subset \eR^n$
\[\Ex \sup_{t\in T} \sum_i t_i g_i g'_i \approx \Ex \sup_{t\in T} \sum_i t_i \exe_i. \]
\end{lemma}
\begin{proof}
Observe that for any $t\geq 1$
\[C^{-1}e^{-Ct} \leq \Pro(|g_i|,|g'_i|\geq \sqrt{t}) \leq \Pro(|g_ig'_i|\geq t) \leq \Pro(|g_i|\geq \sqrt{t})+\Pro(|g_i|\geq \sqrt{t})\leq 2e^{-t/2}. \]
So we can use the same argument as in the proof of Fact \ref{porexp}.
\end{proof}

    \begin{lemma}\label{uff}
    Let $T\subset \eR^n$, $a\in \eR^n$. Define $M=\max_i \sup_{t\in T} |t_i|$. Then
    \[\Ex \sup_{t\in T} \sum_i g_i a_i t^2_i \leq 2M \Ex \sup_{t\in T} \sum_i g_i a_i t_i.\]
        \end{lemma}
        \begin{proof}
            Consider $G_t=\sum_i g_i a_i t^2_i$ and $V_t=2M\sum_i g_i a_i t_i$. Then by trivial calculation we get for any $s,t\in T$ that $\lv G_t-G_{s}\rv^2_2\leq \lv V_t-V_s\rv^2_2$.
            The assertion follows from Slepian's Lemma.
                        \end{proof}

The idea of the next lemma is taken from \cite{adlat} (cf. the proof of Lemma 9.4 therein).

\begin{lemma}\label{kontrak}
Consider $T\subset \eR^{n^2}$ and  $M:=\max_i \sup_{t\in T}\sqrt{\sum_j t^2_{ij}}$. Then
\begin{equation*}
    \Ex \sup_{t\in T} \sum_{ij} g_i g'_i  t_{ij} ^2 \lesssim M \Ex \sup_{t\in T} \sum_{ij} g_{ij}g'_i t_{ij},
\end{equation*}
where $(g_i)_{i \in \eN},(g'_j)_{j\in\eN},(g_{ij})_{i,j\in \eN}$ are independent $\mathcal{N}(0,1)$ r.v.'s.

\end{lemma}
\begin{proof}
By applying Lemma \ref{uff} conditionally on $g'_i$ we get 
\begin{equation} \label{wz1}
\Ex \sup_{t\in T} \sum_{i,j=1}^n g_i g'_i t_{ij}^2\leq 2M\Ex \sup_{t\in T}  \sum_{i=1}^n g_i g'_i \sqrt{ \sum_j t^2_{ij} }.
\end{equation}
We define two  Gaussian processes (conditionally on $g'_1,\ldots,g'_n$) indexed by the set $T$:
\[G_t=\sum_i g_ig'_i \sqrt{\sum_j t^2_{ij}},\qquad V_t=\sum_{i,j} g''_{ij}t_{ij}g'_j.\]
 Cauchy-Schwarz's inequality yields that for any $s,t\in T$, $\Ex^g (G_t-G_s)^2\leq \Ex^{g''} (V_t-V_s)^2$ (conditionally on $g'_1,\ldots,g'_n$). Thus Slepian's Lemma (applied conditionally on $g'_1,\ldots,g'_d$) implies that
\begin{align}
\Ex \sup_{t\in T}  \sum_{i=1}^n g_i g'_i \sqrt{ \sum_j t^2_{ij} }=\Ex^{g'} \Ex^g \sup_{t\in T} G_t \leq \Ex^{g'} \Ex^{g''} \sup_{t\in T} V_t=\Ex \sup_{t\in T} \sum_{ij} g_{ij}g'_j t_{ij}. \label{wz2}
\end{align}
The assertion follows by \eqref{wz1} and \eqref{wz2}. 
\end{proof}

The theorem below was proved in a greater generality by Kwapie{\'n} (for tetrahedral, symmetric polynomials of arbitrary order). We state it in a minimally needed version. For an even more general result (bounds on tails of random chaoses) we refer to \cite{dec1}.
\begin{theorem}\cite[Theorem 2]{dec2}\label{twr:dec}
Let $Q(x,y)=\sum_{i\neq j} a_{ij} x_i y_j$, where $x,y \in \eR^n$ and $a_{ij}=a_{ji}$ are coefficients from  a Banach space $(F,\lv \cdot \rv)$. Let $X_1,\ldots,X_n$ be independent, symmetric r.v.'s. Consider  $Y=(Y_1,\ldots,Y_n)$ an independent copy of $X=(X_1,\ldots,X_n)$. Then for any $p\geq 1$
\[\left( \Ex \lv Q(X,X) \rv^p\right)^{1/p} \approx \left( \Ex \lv Q(X,Y) \rv^p\right)^{1/p}.\]
\end{theorem}

\begin{fact}\label{fak:przek}
Fix $p\geq 1$. Suppose $(a_{ij})_{ij}$ are coefficients from a Banach space $(F,\lv \cdot \rv)$. Let $X_1,\ldots,X_n$ be independent, symmetric r.v.'s such that for any $i$, $\Ex \lv X_i \rv^p <\infty$. Then 
\[\frac{1}{3}\left(\lv \sum_{i\neq j} a_{ij} X_i X_j\rv_p+\lv \sum_{i} a_{i} X^2_i \rv_p \right)\leq \lv \sum_{ij} a_{ij} X_i X_j \rv_p\leq \lv \sum_{i\neq j} a_{ij} X_i X_j\rv_p+\lv \sum_{i} a_{i} X^2_i \rv_p. \]
\end{fact}
\begin{proof}
We prove only the lower bound, the upper bound is obvious. Let $(\eps_i)_i$ be a sequence of i.i.d. Rademacher r.v.'s independent of $(X_i)_i$. By symmetry of $X_i$'s and Jensen's inequality
    \begin{align*}
      &\lv \sum_{ij} a_{ij} X_i X_j \rv_p=\lv \sum_{ij} a_{ij} \eps_i \eps_j X_i X_j \rv_p  \geq \lv \Ex^\eps\sum_{ij} a_{ij} \eps_i \eps_j X_i X_j\rv_p=\lv \sum_{i} a_{i} X^2_i \rv_p.
    \end{align*}
    So by triangle inequality in $L_p$
    \begin{align*}
      \lv \sum_{i\neq j} a_{ij} X_i X_j\rv_p+\lv \sum_{i} a_{i} X^2_i \rv_p\leq   \lv \sum_{ij} a_{ij} X_i X_j\rv_p+2\lv \sum_{i} a_{i} X^2_i \rv_p \leq 3 \lv \sum_{ij} a_{ij} X_i X_j \rv_p.
    \end{align*}
\end{proof}

\section{Glossary}
\begin{center}
	\begin{itemize}
    \item $(F,\lv \cdot \rv)$ --  Banach space with norm $\lv \cdot \rv$,
	\item $\ell^n_q, \ L_q$ -- respectively space of sequences of length $n$ occupied with norm $\lv x \rv^q_q=\sum_i |x_i|^q$ and space of functions integrable to the $q$-th power,
	\item $q,q'$ -- $q$ comes from the underlying $\ell_q$ or $L_q$ space and $q'$ is the H{\"{o}}lder conjugate of $q$ ($1/q+1/q'=1$), 
    \item $(g_i),(g'_i),(g_{ij})$ -- independent $\mathcal{N}(0,1)$ variables,
    \item $G_n$ -- standard normal vector in $\eR^n$,
     \item $(\wy_i),(\wy '_j)$ independent  variables with density $f(x)=1/2 \exp(-|x|/2)$,
    \item $E_n$ -- random vector $(\wy_1,\ldots,\wy_n)$,
    \item $\gamma_{n,\eps},\nu_{n,\eps},\mu_{n\eps}$ -- distributions of respectively $\eps G_n,\ \eps E_n,\ \eps(G_n+E_n)$,
    \item $B_s$ -- set $\{t\in \eR^m : \sum_k |t_k|^s\leq 1\}$,
	\item $N(U,\rho,\eps)$ -- the smallest number of closed balls,
with the diameter $\eps$ in metric $\rho$ that cover the set $U$,
	\item $A$ -- tensor $(a_{ijk})_{i,j\leq n, k\leq m}$,
	\item $v \otimes w$ -- tensor product of vectors $v\in \eR^k$ and $w\in \eR^l$ given by $v \otimes w=(v_i w_j)_{i\leq k, j\leq l} \in \eR^k \times \eR^l$,
	\item $\alpha_A$ -- norm on $\eR^n \times \eR^m$ given by $\alpha_A((x_{jk})_{j\leq n,k\leq m})=\sqrt{\sum_i \left(\sum_{jk} a_{ijk} x_{jk} \right)^2}$,
	\item $d_A$ -- distance on $\eR^n \times \eR^m$ given by 
$d_A((x,t),(x',t'))=\alf(x\otimes t - x' \otimes t')$,
\item $s_A(T)$ -- functional defined on subsets of $\eR^m$  by
\[s_A(T)=\Ex \sup_{t\in T} \left| \sum_{ijk} a_{ijk} g_{ij} t_k\right| + \Ex \sup_{t\in T} \left| \sum_{ijk} a_{ijk} g_i \wy_j t_k \right|, \]
	\item $\beta_{A,T}$ -- norm on $\eR^n$ given by
$\beta_{A,T}(x)=\Ex\sup_{t\in T} \sum_{ijk} a_{ijk} g_i x_j t_k$,
	\item $\Delta_A$ -- diameter with respect to the metric $d_A$,
	\item $F_A(V)$ -- the expected supremum of the Gaussian process indexed by set $V\subset \eR^n \times \eR^m$, given by the formula
\[F_A(V)=\Ex \sup_{(x,t)\in V} \sum_{ijk} a_{ijk} g_i x_j t_k,\]
\item $\alinf$ -- norm on $\eR^n \times \eR^m$ given by 
$\alinf((x_{jk})_{j\leq n,k\leq m})=\max_i \left|\sum_{jk} a_{ijk} x_{jk} \right|$,
	\item $\dinf$ -- distance on $\eR^n \times \eR^m$ given by 
$\dinf((x,t),(x',t'))=\alinf(x\otimes t - x' \otimes t')$,
	\item $\diaminf$ -- diameter with respect do the metric $\dinf$,
	\item $\varphi_A$ -- norm on $\eR^n$ given by $\varphi_A(x)=\sqrt[2q]{\sum_k (\sum_i \frac{(\sum_j a_{ijk}x_j)^4}{\sum_j a^2_{ijk}})^{q/2}}$,
 \item $N^X_i(t), N^Y_i(t)$ -- function given by $N^X_i(t)=-\ln \Pro(|X_i|\geq t)$, formula for $N^Y_j(t)$ is analogous,
	\item r.v.'s with LCT -- class of random variables $X$ such that $t\mapsto -\ln \Pro(|X|\geq t)\in [0,\infty]$ is convex function of $t$,
	\item $\hat{N}^X_i(t), \hat{N}^Y_j(t)$ -- functions given by $\hat{N}^X_i(t)=t^2$ for $|t|\leq 1$ and $\hat{N}^X_i(t)=-\ln \Pro(|X_i|\geq t)$, formula for $\hat{N}^Y_j(t)$ is analogous,
	\item $\norx{(a_i)_i}, \nory{(a_i)_i}$ -- norms on $\eR^n$ given by (the formula for $\nory{(a_i)_i}$ is analogous)
 \[\norx{(a_i)_i}=\sup\{\sum_i a_i x_i : \sum_i \hat{N}^X_i(x_i) \leq p\},\]
\item $\norxy{(a_{ij})_{ij}}$ --  norm on $\eR^n \times \eR^n$ given by 
\[\norxy{(a_{ij})_{ij}}=\sup\{\sum_{ij} a_{ij} x_i y_j : \sum_i \hat{N}^X_i(x_i) \leq p, \sum_j \hat{N}^Y_j(y_j) \leq p \}.\]
    \end{itemize}
\end{center}



\begin{biblist}
\begin{bibdiv}
\bibliographystyle{amsplain}
\bib{Adban}{article}{
   author={Adamczak, Rados\l aw},
   title={Moment inequalities for $U$-statistics},
   journal={Ann. Probab.},
   volume={34},
   date={2006},
   number={6},
   pages={2288--2314},
   review={\MR{2294982}},
}



\bib{adlat}{article}{
   author={Adamczak, Rados{\l}aw},
   author={Lata{\l}a, Rafa{\l}},
   title={Tail and moment estimates for chaoses generated by symmetric
   random variables with logarithmically concave tails},
   journal={Ann. Inst. Henri Poincar\'e Probab. Stat.},
   volume={48},
   date={2012},
   number={4},
   pages={1103--1136},
  review={\MR{3052405}},
}




\bib{adlatmel}{article}{
   author={Adamczak, Rados\l aw},
   author={Lata\l a, Rafa\l },
   author={Meller, Rafa\l },
   title={Hanson-Wright inequality in Banach spaces},
   journal={Ann. Inst. Henri Poincar\'{e} Probab. Stat.},
   volume={56},
   date={2020},
   number={4},
   pages={2356--2376},
   review={\MR{4164840}},
}




\bib{gaban}{article}{
   author={Adamczak, Rados\l aw},
   author={Lata\l a, Rafa\l },
   author={Meller, Rafa\l },
   title={Moments of Gaussian chaoses in Banach spaces},
   journal={Electron. J. Probab.},
   volume={26},
   date={2021},
   pages={Paper No. 11, 36},
  review={\MR{4216524}},
}




\bib{adwolf}{article}{
   author={Adamczak, Rados\l aw},
   author={Wolff, Pawe\l },
   title={Concentration inequalities for non-Lipschitz functions with
   bounded derivatives of higher order},
   journal={Probab. Theory Related Fields},
   volume={162},
   date={2015},
   number={3-4},
   pages={531--586},
   review={\MR{3383337}},
}


\bib{Borell}{collection}{
   author={Christer Borell}
   title={On the Taylor series of a Wiener polynomial,},
   series={Seminar Notes on multiple stochastic integration,
polynomial chaos and their integration},
publisher={Case Western Reserve Univ., Cleveland}
   date={1984},
   review={\MR{4164840}},
}

\bib{dec1}{article}{
   author={de la Pe{\~n}a, Victor H.},
   author={Stephen Montgomery-Smith},
   title={Decoupling inequalities for the tail probabilities of multivariate
   $U$-statistics},
   journal={Ann. Probab.},
   volume={23},
   date={1995},
   number={2},
   pages={806--816},
   review={\MR{1334173}},
}

\bib{nayar}{article}{
   author={Eskenazis, Alexandros},
   author={Nayar, Piotr},
   author={Tkocz, Tomasz},
   title={Gaussian mixtures: entropy and geometric inequalities},
   journal={Ann. Probab.},
   volume={46},
   date={2018},
   number={5},
   pages={2908--2945},
  review={\MR{3846841}},
}
	



\bib{glu}{article}{
   author={Efim Gluskin},
   author={Stanis{\l{}}aw Kwapie{\'n}},
   title={Tail and moment estimates for sums of independent random variables
   with logarithmically concave tails},
   journal={Studia Math.},
   volume={114},
   date={1995},
   number={3},
   pages={303--309},
   review={\MR{1338834}},
}
\bib{Gotze}{article}{
   author={G\"{o}tze, Friedrich},
   author={Sambale, Holger},
   author={Sinulis, Arthur},
   title={Concentration inequalities for polynomials in
   $\alpha$-sub-exponential random variables},
   journal={Electron. J. Probab.},
   volume={26},
   date={2021},
   pages={Paper No. 48, 22},
   review={\MR{4247973}},
}



\bib{Handbook}{collection}{
   title={Handbook of the geometry of Banach spaces. Vol. I},
   editor={William Buhmann Johnson},
   editor={Joram Lindenstrauss},
   publisher={North-Holland Publishing Co., Amsterdam},
   date={2001},
   review={\MR{1863688}},
}


\bib{kol}{article}{
   author={Kolesko, Konrad},
   author={Lata\l a, Rafa\l },
   title={Moment estimates for chaoses generated by symmetric random
   variables with logarithmically convex tails},
   journal={Statist. Probab. Lett.},
   volume={107},
   date={2015},
   pages={210--214},
   review={\MR{3412778}},
}


\bib{dec2}{article}{
   author={Kwapie{\'n}, Stanis{\l}aw},
   title={Decoupling inequalities for polynomial chaos},
   journal={Ann. Probab.},
   volume={15},
   date={1987},
   number={3},
   pages={1062--1071},
   review={\MR{893914}},
}




\bib{latgaus}{article}{
   author={Lata{\l}a, Rafa{\l}},
   title={Estimates of moments and tails of Gaussian chaoses},
   journal={Ann. Probab.},
   volume={34},
   date={2006},
   number={6},
   pages={2315--2331},
   review={\MR{2294983}},
}


\bib{d1}{article}{
   author={Lata\l a, Rafa\l },
   title={Estimation of moments of sums of independent real random
   variables},
   journal={Ann. Probab.},
   volume={25},
   date={1997},
   number={3},
   pages={1502--1513},
   review={\MR{1457628}},
}


\bib{latsud}{article}{
   author={Lata{\l}a, Rafa{\l}},
   title={Sudakov minoration principle and supremum of some processes},
   journal={Geom. Funct. Anal.},
   volume={7},
   date={1997},
   number={5},
   pages={936--953},
   review={\MR{1475551}},
}



\bib{some}{article}{
   author={Lata\l a, Rafa\l },
   title={Tail and moment estimates for some types of chaos},
   journal={Studia Math.},
   volume={135},
   date={1999},
   number={1},
   pages={39--53},
   review={\MR{1686370}},
}

\bib{lat}{article}{
   author={Lata{\l}a, Rafa{\l}},
   title={Tail and moment estimates for sums of independent random vectors
   with logarithmically concave tails},
   journal={Studia Math.},
   volume={118},
   date={1996},
   number={3},
   pages={301--304},
   review={\MR{1388035}},
}
\bib{latroy}{article}{
   author={Lata\l a, Rafa\l },
   author={Matlak, Dariusz},
   title={Royen's proof of the Gaussian correlation inequality},
   conference={
      title={Geometric aspects of functional analysis},
   },
   book={
      series={Lecture Notes in Math.},
      volume={2169},
      publisher={Springer, Cham},
   },
   date={2017},
   pages={265--275},
   review={\MR{3645127}},
}

\bib{strzelec}{article}{
   author={Lata\l a, Rafa\l },
   author={Strzelecka, Marta},
   title={Comparison of weak and strong moments for vectors with independent
   coordinates},
   journal={Mathematika},
   volume={64},
   date={2018},
   number={1},
   pages={211--229},
   review={\MR{3778221}},
}

\bib{lattkocz}{article}{
   author={Lata\l a, Rafa\l },
   author={Tkocz, Tomasz},
   title={A note on suprema of canonical processes based on random variables
   with regular moments},
   journal={Electron. J. Probab.},
   volume={20},
   date={2015},
   pages={no. 36, 17},
  review={\MR{3335827}},
}


\bib{proinban}{book}{
   author={Ledoux, Michel},
   author={Talagrand, Michel},
   title={Probability in Banach spaces},
   series={Ergebnisse der Mathematik und ihrer Grenzgebiete (3) [Results in
   Mathematics and Related Areas (3)]},
   volume={23},
   publisher={Springer-Verlag, Berlin},
   date={1991},
   review={\MR{1102015}},
}

\bib{loch}{article}{
   author={Rafa\l \  Marcin  \L ochowski },
   title={Moment and tail estimates for multidimensional chaoses generated
   by symmetric random variables with logarithmically concave tails},
   conference={
      title={Approximation and probability},
   },
   book={
      series={Banach Center Publ.},
      volume={72},
      publisher={Polish Acad. Sci. Inst. Math., Warsaw},
   },
   date={2006},
   pages={161--176},
   review={\MR{2325744}},
}


\bib{MR3984282}{article}{
   author={Meller, Rafa\l },
   title={Tail and moment estimates for a class of random chaoses of order
   two},
   journal={Studia Math.},
   volume={249},
   date={2019},
   number={1},
   pages={1--32},
   review={\MR{3984282}},
}

\bib{ja}{article}{
   author={Meller, Rafa\l },
   title={Two-sided moment estimates for a class of nonnegative chaoses},
   journal={Statist. Probab. Lett.},
   volume={119},
   date={2016},
   pages={213--219},
  review={\MR{3555289}},
}




\bib{royen}{article}{
   author={Royen, Thomas},
   title={A simple proof of the Gaussian correlation conjecture extended to
   some multivariate gamma distributions},
   journal={Far East J. Theor. Stat.},
   volume={48},
   date={2014},
   number={2},
   pages={139--145},
   review={\MR{3289621}},
}
\bib{tal}{book}{
   author={Talagrand, Michel},
   title={Upper and lower bounds for stochastic processes---decomposition
   theorems},
   series={Ergebnisse der Mathematik und ihrer Grenzgebiete. 3. Folge. A
   Series of Modern Surveys in Mathematics [Results in Mathematics and
   Related Areas. 3rd Series. A Series of Modern Surveys in Mathematics]},
   volume={60},
   note={Second edition [of  3184689]},
   publisher={Springer, Cham},
   date={[2021] \copyright 2021},
   review={\MR{4381414}},
}
	



\end{bibdiv}

\end{biblist}
\end{document}